\documentclass[12pt]{amsart}
\usepackage{amssymb}
\usepackage{amsmath}

\theoremstyle{plain}
\newtheorem{theorem}{Theorem}[section]
\newtheorem{corollary}[theorem]{Corollary}

\newtheorem{lemma}[theorem]{Lemma}
\newtheorem{proposition}[theorem]{Proposition}

\theoremstyle{definition}

\newtheorem{example}[theorem]{Example}

\theoremstyle{remark}
\newtheorem{remark}[theorem]{Remark}

\def\Z{{\mathbb Z}}
\def\N{{\mathbb N}}
\def\R{{\mathbb R}}
\def\C{{\mathbb C}}
\def\Q{{\mathbb Q}}
\def\I{{\mathbb I}}
\def\V{{\mathbb V}}

\def\cF{{\mathcal F}}
\def\cO{{\mathcal O}}

\def\fkm{{\mathfrak m}}

\def\a{{\boldsymbol a}}
\def\b{{\boldsymbol b}}
\def\c{{\boldsymbol c}}

\def\l{{\boldsymbol l}}
\def\m{{\boldsymbol m}}

\def\u{{\boldsymbol u}}
\def\v{{\boldsymbol v}}
\def\x{{\boldsymbol x}}
\def\y{{\boldsymbol y}}

\def\0{{\boldsymbol 0}}
\def\1{{\boldsymbol 1}}

\def\aalpha{{\boldsymbol{\alpha}}}
\def\bbeta{{\boldsymbol{\beta}}}
\def\ggamma{{\boldsymbol{\gamma}}}

\def\llambda{{\boldsymbol{\lambda}}}
\def\mmu{{\boldsymbol{\mu}}}

\def\End{{\rm End}}
\def\Hom{{\rm Hom}}
\def\Ann{{\rm Ann}}

\def\ZC{{\rm ZC}}
\def\Prim{{\rm Prim}}

\title[Primitive ideals of the ring of differential operators]{Primitive 
ideals of the ring of differential operators
on an affine toric variety}
\author{Mutsumi Saito}

\date{May 31, 2005}

\begin{document}


\begin{abstract}
Let $A$ be a $d\times n$ integer matrix whose column vectors
generate the lattice $\Z^d$,
and let $D(R_A)$ be the ring of differential operators on the
affine toric variety defined by $A$.

We show that the classification of $A$-hypergeometric systems
and that of $\Z^d$-graded simple $D(R_A)$-modules (up to shift) are the same.
We then show that
the set of $\Z^d$-homogeneous primitive ideals of $D(R_A)$ is finite.
Furthermore, we give conditions for the algebra $D(R_A)$ being simple.


\end{abstract}

\maketitle

\section{Introduction}

Let $A$ be a $d\times n$ integer matrix whose column vectors
generate the lattice $\Z^d$.
Let $R_A$ be the ring of regular functions
on the
affine toric variety defined by $A$, and
$D(R_A)$ its ring of differential operators.

In this paper, we show the following three theorems:
\begin{enumerate}
\item
The classification of 
$A$-hypergeometric systems
and that of $\Z^d$-graded simple $D(R_A)$-modules (up to shift)
are the same (Theorem \ref{thm:summary}).
\item
The set of $\Z^d$-homogeneous primitive ideals of $D(R_A)$ is finite
(Theorem \ref{thm:PrimFinite}).
\item
The algebra $D(R_A)$ is simple if and only if
$R_A$ is a scored semigroup ring, and $A$ satisfies
a certain condition (C2) 
(Theorem \ref{simple<=>scored,C2}).
\end{enumerate}

The ring of differential operators was introduced by 
Grothendieck \cite{ega4} and Sweedler \cite{Sw}. 
As for the ring of differential operators $D(R_A)$ 
on an affine toric variety,
many recent papers such as
Jones \cite{Jones}, Musson \cite{Musson}, and Musson and Van
den Bergh \cite{Musson-Van den Bergh}
describe the structure of $D(R_A)$ when $R_A$ is normal.
For general $R_A$, we studied the finite generation of
$D(R_A)$ and its graded ring with respect to the order filtration
in \cite{Saito-Traves} and \cite{Saito-Traves2}.
Moreover, in \cite{Saito-Traves}, we showed that the algebra $D(R_A)$ 
and the symmetry algebra (the algebra of contiguity operators) 
of $A$-hypergeometric systems are anti-isomorphic to each other.
This paper may be considered as a continuation of \cite{Saito-Traves}.

The history of $A$-hypergeometric systems (or GKZ hypergeometric
systems) goes back to
Kalnins, Manocha, and Miller \cite{KMM}, and Hrabowski \cite{Hrabowski}.
After the papers by
Gel'fand, Kapranov, and Zelevinskii (e.g. \cite{Gelfand}--\cite{GKZ}), 
researchers in various fields have studied the systems,
and established connection with representation theory,
algebraic geometry, commutative ring theory, etc.
See, for example, the bibliography of \cite{SST}.

Associated to a parameter vector $\aalpha$ 
and a face $\tau$ of the cone generated by
the column vectors of $A$, we defined a finite set $E_\tau(\aalpha)$
(see \eqref{def:Etau} in Section 7) in
\cite{IsoClass}, and proved that the $A$-hypergeometric systems are
classified by those finite sets.
Hence in order to show that the classification of 
$A$-hypergeometric systems
and that of $\Z^d$-graded simple $D(R_A)$-modules (up to shift) are the same,
we only need to show that $\Z^d$-graded simple $D(R_A)$-modules
can be classified (up to shift) by the finite sets $E_\tau(\aalpha)$
as shown in Theorem \ref{thm:Liffsim},
essentially in \cite{Musson-Van den Bergh} and \cite{Saito-Traves}.
It is however desirable to show the equivalence
intrinsically.
We therefore present the second proof of the equivalence
by connecting $A$-hypergeometric systems with certain
$\Z^d$-graded $D(R_A)$-modules by functors
(Corollary \ref{cor:MRiffHR}).

Some topics discussed in this paper were treated 
in \cite{Musson-Van den Bergh} under the conditions $(A1)$ and $(A2)$
(see p.4 in \cite{Musson-Van den Bergh}).
In our case, $(A1)$ is always satisfied, and
$(A2)$ requiring that all $\Z^d$-homogeneous components $D(R_A)_\a$ are 
singly generated $D(R_A)_\0$-modules
is equivalent to requiring
that $R_A$ satisfies the Serre's $(S_2)$ condition
(see Proposition \ref{A2=S2}).
In this paper, we do not assume the Serre's $(S_2)$ condition.

The layout of this paper is as follows:
we start with recalling the definitions and some fundamental facts
about differential operators in Section 2, and about $A$-hypergeometric
systems in Section 3.
In Section 4, we recall some results
by Musson and Van den Bergh \cite{Musson-Van den Bergh}
about the category $\cO$, analogous to the Bernstein-Gel'fand-Gel'fand's
category $\cO$, and the counterpart ${}^R\cO$ for right $D(R_A)$-modules.
Then we recall realizations $L(\aalpha)$ (${}^RL(\aalpha)$)
of simple objects in $\cO$ (${}^R\cO$)
from \cite{Saito-Traves}, and we
show that we have a duality between $\cO$ and ${}^R\cO$
sending $L(\aalpha)$ and ${}^RL(\aalpha)$ to each other.
We also show that any $\Z^d$-graded simple $D(R_A)$-module is
isomorphic (up to shift) to $L(\aalpha)$.
These combined prove that the classifications of $L(\aalpha)$,
${}^R L(\aalpha)$, their projective covers $M(\aalpha)$,
${}^R M(\aalpha)$, and $\aalpha$ with respect to the equivalence
relation determined by the finite sets $E_\tau(\aalpha)$
are the same.

In Section 5, we provide two functors between the category of
right $D(R_A)$-modules and the category of right $D(R)$-modules
supported by the affine toric variety defined by $A$, where
$D(R)$ is the $n$-th Weyl algebra.
One is the direct image functor, for right $D$-modules,
of the closed inclusion of the affine toric variety into $\C^n$,
and the other is its right adjoint functor.
By using these functors,
we prove that ${}^R M(\aalpha)\simeq {}^R M(\bbeta)$
if and only if their corresponding $A$-hypergeometric systems 
are isomorphic (Theorem \ref{thm:summary}).

After we see a couple of basic facts about primitive ideals
in Section 6,
we show in Section 7 that if we properly perturb a parameter $\aalpha$
then the annihilator ideal $\Ann\, L(\aalpha)$ remains unchanged,
and in this way we show that the set $\Prim\, D(R_A)$ of 
$\Z^d$-homogeneous primitive ideals 
of $D(R_A)$ is finite (Theorem \ref{thm:PrimFinite}).

In Section 8, the simplicity of $D(R_A)$ is treated.
First we consider the conditions: the scoredness and the Serre's $(S_2)$.
We prove that the simplicity of $D(R_A)$ implies the scoredness, and that
the conditions $(A2)$ and $(S_2)$ are equivalent.
Finally we give a necessary and sufficient condition for the simplicity
(Theorem \ref{simple<=>scored,C2}).
In Section 9, we give an example of computation of 
the set $\Prim\, D(R_A)$, and an example
such that $R_A$ is scored and Cohen-Macaulay, but $D(R_A)$ is not simple.

\section{Ring of differential operators
on an affine toric variety}

In this section, 
we recall some fundamental facts about
the rings of differential operators of semigroup algebras.

Let
$
A:=\{\, {\a}_1, {\a}_2,\ldots,
{\a}_n\,\}
$
be a finite set of column vectors in $\mathbb{Z}^d$.  Sometimes we identify
$A$ with the matrix $(\a_1,\a_2,\ldots,\a_n)$.
Let $\N A$ and $\Z A$ denote the monoid and the abelian group
generated by $A$, respectively.
Throughout this paper, we assume that $\mathbb{Z} A=\mathbb{Z}^d$ for simplicity.

Let $R$ denote the polynomial ring $\C[x]:=\C[x_1,\ldots, x_n]$.
The semigroup algebra $R_A:=\mathbb{C}[\mathbb{N} A]=\bigoplus_{\a\in \mathbb{N} A}\mathbb{C} t^\a$ is
the ring of regular functions on the affine toric variety defined by
$A$, where $t^\a=t_1^{a_1}t_2^{a_2}\cdots t_d^{a_d}$ for
$\a={}^t(a_1,a_2,\ldots,a_d)$. 
Then we have $R_A\simeq R/I_A(x)$, where
$I_A(x)$ is the ideal of $\C[x]$
generated by all $x^\u-x^\v$
for $\u, \v\in \N^n$ with $A\u=A\v$.

Let $M, N$ be $R$-modules.
We briefly recall the module $D(M,N)$ of differential operators
from $M$ to $N$. For details, see \cite{Smith-Stafford}.
For $k\in \N$, 
the subspaces $D^k(M,N)$ of $\Hom_\C(M,N)$
are inductively defined by
$$ D^0(M,N)=\Hom_R(M,N)$$ and
$$
 D^{k+1}(M,N)=\{ P\in \Hom_\C(M,N)\, :\,
[f, P]\in D^k(M,N)\quad (\forall f\in R\},
$$
where $[\,,\,]$ denotes the commutator.
Set $D(M,N):=\bigcup_{k=0}^\infty D^k(M,N)$,
and $D(M):=D(M,M)$.
Then $D(M)$ is a $\C$-algebra, and
$D(M,N)$ is a $(D(N), D(M))$-bimodule.
Hence $D(R, R_A)$ is a $(D(R_A), D(R))$-module.

The ring $D(R)$ is the $n$-th Weyl algebra
$$
D(R)=\C \langle x_1,\ldots, x_n,
\frac{\partial}{\partial x_1},\ldots,\frac{\partial}{\partial x_n} 
\rangle,
$$
where $[\frac{\partial}{\partial x_i},x_j]=\delta_{ij}$,
and the other
pairs of generators commute.
Here $\delta_{ij}$ is 1 if $i=j$ and 0 otherwise.

Let $\C[t, t^{-1}]$ denote the Laurent polynomial
ring $\C[t_1^{\pm 1},\ldots, t_d^{\pm 1}]$.
Then its ring of differential operators $D(\C[t, t^{-1}])$ is
the ring 
$$
\mathbb{C}[t^{\pm 1}_1,\ldots,t^{\pm 1}_d]\langle \partial_1,
\ldots,\partial_d\rangle,
$$
where $[\partial_i, t_j]=\delta_{ij}$, 
$[\partial_i, t_j^{-1}]=-\delta_{ij}t_j^{-2}$,
and the other
pairs of generators commute.
The ring of differential operators $D(R_A)$ can be
realized as a subring of the ring $D(\C[t, t^{-1}])$
by
$$
D(R_A)=\{
P\in 
\mathbb{C}[t^{\pm 1}_1,\ldots,t^{\pm 1}_d]\langle \partial_1,
\ldots,\partial_d\rangle
\, :\, P(R_A)\subseteq R_A\}.
$$

Put $s_j:=t_j\partial_j$ for $j=1,2,\ldots, d$.
Then it is easy to see that $s_j\in D(R_A)$ for all $j$.
We introduce a $\Z^d$-grading on the ring $D(R_A)$;
for $\a={}^t(a_1,a_2,\ldots,a_d)\in \Z^d$, set
$$
D(R_A)_\a:=
\{ P\in D(R_A)\, :\, [s_j, P]=a_jP\quad\mbox{for $j=1,2,\ldots,d$}\}.
$$
Then $D(R_A)=\bigoplus_{\a\in \Z^d}D(R_A)_\a$.

By regarding $D(R_A)_\a$ as a subset of 
$\mathbb{C}[t^{\pm 1}_1,\ldots,t^{\pm 1}_d]\langle \partial_1,
\ldots,\partial_d\rangle$,
we see that there exists an ideal $I$ of $\C[s]:=\C[s_1,\ldots,s_d]$
such that $D(R_A)_\a=t^\a I$.
To describe this ideal $I$ explicitly, 
we define a subset
$\Omega(\a)$ of the semigroup $\N A$ by
$$
\Omega(\a) = \{\, \b
\in \N A: \; \b+ \a \not\in \N A\,\}
=
\N A\setminus (-\a +\N A).
$$
Then each $D(R_A)_\a$ is described as follows.

\begin{theorem}[\cite{Jones}, Theorem 3.3.1 in \cite{Saito-Traves}]
\label{thm:Jones}
$$D(R_A)_\a=
 t^{\a} {\I} (\Omega(\a)) \mbox{\quad for all $\a\in \Z^d$,}
$$
where
$$
\I(\Omega(\a)):=\{
f(s)\in \C[s]\,:\,
\mbox{$f$ vanishes on $\Omega(\a)$}\}.
$$ \label{thm3.3.1}
\end{theorem}

In particular, $D(R_A)_\a=t^\a \C[s]=\C[s] t^\a$
for each $\a\in \N A$,
since $\Omega(\a)=\emptyset$ in this case.

\section{$A$-hypergeometric systems}

Let us briefly recall the definition of an $A$-hypergeometric system
and its classification.

Let $\aalpha={}^t(\alpha_1,\ldots,\alpha_d)\in \C^d$.
The $A$-hypergeometric system with parameter $\aalpha$ is
the left $D(R)$-module
$$
H_A(\aalpha)
:=D(R)/\left(
\sum_{i=1}^d D(R)(\sum_{j=1}^n a_{ij}x_j\frac{\partial}{\partial x_j}
-\alpha_i)+D(R)I_A(\partial)\right),
$$
where $\a_j={}^t(a_{1j}, a_{2j},\ldots, a_{dj})$,
$I_A(\partial)$ is the ideal of $\C[\frac{\partial}{\partial x_1},\ldots,\frac{\partial}{\partial x_n}]$
generated by all $\prod_{j=1}^n \frac{\partial^{u_j}}{\partial^{u_j} x_j}
-\prod_{j=1}^n \frac{\partial^{v_j}}{\partial^{v_j} x_j}$
for $\u, \v\in \N^n$ with $A\u=A\v$.

Interchanging $x_j$ and $\frac{\partial}{\partial x_j}$ for all $j$,
we have an anti-automorphism $\iota$ of $D(R)$.
Clearly $\iota$ gives rise to a one-to-one correspondence between
the left $D(R)$-modules and the right $D(R)$-modules.
Thus $\iota$ induces a right $D(R)$-module
$$
{}^R H_A(\aalpha)
:=D(R)/\left(
\sum_{i=1}^d (\sum_{j=1}^n a_{ij}x_j\frac{\partial}{\partial x_j}
-\alpha_i)D(R)+I_A(x) D(R)\right).
$$
Note that $\iota(x_j\frac{\partial}{\partial x_j})=
\iota(\frac{\partial}{\partial x_j})\iota(x_j)=
x_j\frac{\partial}{\partial x_j}$.

In \cite[Definition 4.1.1]{Saito-Traves},
we have introduced a partial order into the parameter space $\C^d$
(see \eqref{def:Preceq}), 
which is equivalent
by \cite[Lemma 4.1.4]{Saito-Traves} to
\begin{equation}
\label{eq:PartialOrder1}
\aalpha\preceq \bbeta
\Longleftrightarrow
\I(\Omega(\bbeta-\aalpha))\not\subseteq \fkm_\aalpha,
\end{equation}
where $\fkm_\aalpha$ is the maximal ideal of $\C[s]$ at $\aalpha$.
Note that, if $\bbeta-\aalpha\notin \Z^d$, then $\Omega(\bbeta-\aalpha)
=\N A$, and hence $\aalpha \not\preceq \bbeta$.
We write $\aalpha\sim\bbeta$ if $\aalpha\preceq\bbeta$ and
$\aalpha\succeq\bbeta$.
This equivalence relation was introduced also by
Musson and Van den Bergh 
(see \cite[Lemma 3.1.9 (6)]{Musson-Van den Bergh}).

This relation classifies $A$-hypergeometric systems.

\begin{theorem}[Theorem 2.1 in \cite{IsoClass}]
\label{thm:ClassificationOfAHGS}
$H_A(\aalpha)\simeq H_A(\bbeta)$
if and only if $\aalpha\sim\bbeta$.
\end{theorem}

\section{$D(R_A)$-modules}

In this section, we recall some results
by Musson and Van den Bergh \cite{Musson-Van den Bergh}
about the category $\cO$, 
and the counterpart ${}^R\cO$ for right $D(R_A)$-modules.
Then we recall realizations $L(\aalpha)$ (${}^RL(\aalpha)$)
of simple objects in $\cO$ (${}^R\cO$)
from \cite{Saito-Traves}, and we
show that we have a duality between $\cO$ and ${}^R\cO$
sending $L(\aalpha)$ and ${}^RL(\aalpha)$ to each other.
We also show that any $\Z^d$-graded simple $D(R_A)$-module is
isomorphic (up to shift) to $L(\aalpha)$.
These combined prove that the classifications of $L(\aalpha)$,
${}^R L(\aalpha)$, their projective covers $M(\aalpha)$,
${}^R M(\aalpha)$, and $\aalpha$ with respect to the equivalence
relation $\sim$
are the same.

\subsection{Left modules}

Let us recall the full subcategory $\cO$ of the category
of left $D(R_A)$-modules introduced in \cite{Musson-Van den Bergh},
which is an analogue of the Bernstein-Gel'fand-Gel'fand's category
$\cO$ for the study of highest weight modules of semisimple Lie
algebras.
A left $D(R_A)$-module $M$ is an object of $\cO$ if
$M$ has a weight decomposition
$M=\bigoplus_{\llambda\in \C^d}M_\llambda$
with each $M_\llambda$ finite-dimensional,
where
$$
M_\llambda=\{
x\in M\, :\, f(s). x= f(\llambda)x\quad \text{(for all $f\in \C[s]$)}\}.
$$
We call $\llambda$ a weight of $M$ if $M_\llambda\neq 0$.

For $\aalpha={}^t(\alpha_1,\ldots,\alpha_d)\in \C^d$,
set 
$$
M(\aalpha):=D(R_A)/D(R_A)(s-\aalpha),
$$
where $D(R_A)(s-\aalpha)$ means
$\sum_{i=1}^d D(R_A)(s_i-\alpha_i)$.
Then
$M(\aalpha)=\bigoplus_{\llambda\in \aalpha+\Z A}M(\aalpha)_\llambda$,
and
$M(\aalpha)\in \cO$.

Among others, Musson and Van den Bergh proved the following.

\begin{proposition}[Proposition 3.1.7 in \cite{Musson-Van den Bergh}]
\label{prop:Musson-VandenBergh}
\qquad\qquad\qquad\qquad\qquad
\begin{enumerate}
\item
${\rm Hom}_{D(R_A)}(M(\aalpha), M)=M_\aalpha$
for $M\in \cO$.
\item
$M(\aalpha)$ is a projective object in $\cO$.
\item
$M(\aalpha)$ has a unique simple quotient
(denoted by $L(\aalpha)$).
\item
All simple objects in $\cO$ are of the form $L(\aalpha)$.
\item
The natural projection $M(\aalpha)\to L(\aalpha)$ is the projective cover.
\item
$M(\aalpha)\simeq M(\bbeta)$
if and only if
$L(\aalpha)\simeq L(\bbeta)$.
\end{enumerate}
\end{proposition}

\begin{remark}
\label{rem:A1andA2}
Musson and Van den Bergh
assumed the conditions $(A1)$ and $(A2)$
(see p.4 in \cite{Musson-Van den Bergh}).
In our case, $(A1)$ is always satisfied, and
$(A2)$ requiring that all $D(R_A)_\a$ are 
singly generated $\C [s]$-modules
is equivalent to requiring
that $R_A$ satisfies the Serre's $(S_2)$ condition
(see Proposition \ref{A2=S2}).
For Proposition \ref{prop:Musson-VandenBergh}, we do not need the condition $(A2)$.
\end{remark}

\begin{remark}
\label{rem:SimpleOmeansSimple}
Let $M\in \cO$, and let
$N$ be a left $D(R_A)$-submodule of $M$.
Then $N\in \cO$.
Hence, if $M$ is a simple object in $\cO$,
then $M$ is a simple left $D(R_A)$-module.
\end{remark}

Let $\aalpha\in \C^d$.
In \cite{Saito-Traves}, we studied the composition factors of a $D(R_A)$-module
$$
\C[t_1^{\pm 1},\ldots, t_d^{\pm 1}]t^\aalpha,
$$
and
we saw that
\begin{equation}
\label{eq:LeftL(alpha)}
\bigoplus_{\llambda\succeq\aalpha}\C t^\llambda/
\bigoplus_{\llambda\succ\aalpha}\C t^\llambda
\end{equation}
is simple \cite[Theorem 4.1.6]{Saito-Traves},
where $\llambda\succ\aalpha$ means $\llambda\succeq\aalpha$
and $\llambda\not\sim\aalpha$.
The $D(R_A)$-module \eqref{eq:LeftL(alpha)} 
is a simple quotient of $M(\aalpha)$,
and hence a realization of $L(\aalpha)$.
In particular, the set of weights of $L(\aalpha)$ is
$$
\{ \llambda\in \C^d\,:\, \llambda\sim \aalpha\}.
$$

\begin{theorem}[cf. Lemma 3.1.9 (6) in \cite{Musson-Van den Bergh}]
\label{thm:Liffsim}
$L(\aalpha)\simeq L(\bbeta)$ if and only if $\aalpha\sim\bbeta$.
\end{theorem}

\begin{proof}
By the realization \eqref{eq:LeftL(alpha)},
$L(\aalpha)= L(\bbeta)$ if $\aalpha\sim\bbeta$.

If $\aalpha\not\sim\bbeta$, then $L(\aalpha)$ and $L(\bbeta)$
have different weights. Hence $L(\aalpha)\not\simeq L(\bbeta)$.
\end{proof}

\begin{proposition}
Let $M$ be a $\Z^d$-graded simple left $D(R_A)$-module.
Then $M$ is isomorphic to $L(\aalpha)$ for some $\aalpha$
as a left $D(R_A)$-module.
\end{proposition}

\begin{proof}
Recall that $D(R_A)_\0=\C [s]$.
First we show that
each nonzero $M_\llambda$ is a simple $\C[s]$-module.
Suppose that $N_\llambda$ is a nontrivial $\C[s]$-submodule of $M_\llambda$. Put $N:=\bigoplus_{\a\in \Z^d}N_{\llambda+\a}=
\bigoplus_{\a\in \Z A}D(R_A)_{\a}N_{\llambda}$.
Then $N$ is a nontrivial $\Z^d$-graded submodule of $M$,
which contradicts the assumption.

Suppose that $M_\llambda\neq 0$.
By the first paragraph, there exists $\aalpha\in \C^d$ such that
$M_\llambda\simeq \C[s]/\fkm_\aalpha$ as $\C[s]$-modules.
Let $M[\llambda-\aalpha]$ be the $\Z^d$-graded $D(R_A)$-module
shifted by $\llambda-\aalpha$, i.e.,
$M[\llambda-\aalpha]_\mmu=M_{\mmu+\llambda-\aalpha}$.
Then $M[\llambda-\aalpha]=\bigoplus_{\mmu\in \aalpha+\Z^d}M[\llambda-\aalpha]_\mmu$, and 
$M[\llambda-\aalpha]_{\aalpha+\a}=D(R_A)_\a M_\llambda$.
Hence $M[\llambda-\aalpha]\in \cO$, and
$M[\llambda-\aalpha]\simeq L(\aalpha)\in\cO$.
\end{proof}

\subsection{Right modules}

A right $D(R_A)$-module $M$ is an object of 
${}^R\cO$ if
$M$ has a weight decomposition
$M=\bigoplus_{\llambda\in \C^d}M_\llambda$
with each $M_\llambda$ finite-dimensional,
where
$$
M_\llambda=\{
x\in M\, :\, x.f(s)= f(-\llambda)x\quad (\forall f\in \C[s])\}.
$$
We can make a parallel argument
about the categories $\cO$ and ${}^R\cO$.
Indeed we shall show that there exists a duality functor between them.

For $\aalpha\in \C^d$,
set 
$$
{}^R M(\aalpha):=D(R_A)/(s-\aalpha)D(R_A).
$$
Then
${}^R M(\aalpha)=\bigoplus_{\llambda\in -\aalpha+\Z A}
{}^R M(\aalpha)_\llambda$,
and
${}^R M(\aalpha)\in {}^R\cO$.

The following proposition is proved in the same way 
as Proposition \ref{prop:Musson-VandenBergh}.

\begin{proposition}
\label{prop:Musson-VandenBerghR}

\begin{enumerate}
\item
${\rm Hom}_{D(R_A)}({}^R M(\aalpha), M)=M_{-\aalpha}$
for $M\in {}^R\cO$.
\item
${}^R M(\aalpha)$ is a projective object in ${}^R \cO$.
\item
${}^R M(\aalpha)$ has a unique simple quotient
(denoted by ${}^R L(\aalpha)$).
\item
All simple objects in ${}^R\cO$ are of the form ${}^R L(\aalpha)$.
\item
The natural projection ${}^R M(\aalpha)\to {}^R L(\aalpha)$ is the projective cover.
\item
${}^R M(\aalpha)\simeq {}^R M(\bbeta)$
if and only if
${}^R L(\aalpha)\simeq {}^R L(\bbeta)$.
\end{enumerate}
\end{proposition}

The ring $D(R_A)$ is a subring of 
$\C[t_1^{\pm 1},\ldots, t_d^{\pm 1}]\langle
\partial_1, \ldots, \partial_d \rangle$,
where we can take formal adjoint operators, and thus
we can consider a right $D(R_A)$-module
$$
\C[t_1^{\pm 1},\ldots, t_d^{\pm 1}]t^\aalpha \frac{dt}{t}.
$$
Here the right action of $P=\sum_\a t^\a f_\a(s)$
on this module is defined by
$$
\left( g(t)\frac{dt}{t}\right).P:=
P^*(g)\frac{dt}{t},
$$
where
$P^*=\sum_\a f_\a(-s)t^\a$,
and recall that $s_i=t_i\partial_i$ ($i=1,\ldots d$).

Then
\begin{equation}
\bigoplus_{\bbeta\preceq\aalpha}\C t^{-\bbeta}\frac{dt}{t}
/
\bigoplus_{\bbeta\prec\aalpha}\C t^{-\bbeta}\frac{dt}{t}.
\end{equation}
is a realization of ${}^R L(\aalpha)$.

Let $M\in \cO$ $({}^R\cO)$.
Then $\Hom_\C( M,\C)$ is a right (left) $D(R_A)$-module.
Define a right (left) $D(R_A)$-submodule $M^*$ of
$\Hom_\C( M,\C)$ by
$$
M^*:=\bigoplus_{\llambda}M^*_\llambda,\quad
M^*_\llambda:=
\Hom_\C( M_{-\llambda}, \C).
$$
Then ${}^*:\cO\to {}^R\cO$ $({}^*:{}^R\cO\to \cO)$ is a duality
functor.
Hence we have the following proposition.

\begin{proposition}
\label{prop:duality}
\begin{enumerate}
\item
${\rm Hom}_{D(R_A)}(M, {}^R M(\aalpha)^*)=
\Hom_\C( M_{\aalpha}, \C)$
for $M\in \cO$.
\item
${}^R M(\aalpha)^*$ is an injective object in $\cO$.
\item
${}^R M(\aalpha)^*$ has a unique simple subobject 
${}^R L(\aalpha)^*$in $\cO$.
\item
$L(\aalpha)\simeq {}^R L(\aalpha)^*$.
\item
The natural inclusion ${}^R L(\aalpha)^*\to {}^R M(\aalpha)^*$ is the
injective hull.
\item
$L(\aalpha)\simeq L(\bbeta)$
if and only if
${}^R L(\aalpha)\simeq {}^R L(\bbeta)$.
\end{enumerate}
\end{proposition}

\begin{proof}
(4) follows from the fact that two simple modules
$L(\aalpha)$ and ${}^R L(\aalpha)^*$
have the same weight spaces.

The other statements are clear.
\end{proof}

\section{$A$-hypergeometric systems and Category $\mathcal{O}$}
\label{Section:CategoryO}

We have proved that 
the classification of $A$-hypergeometric
systems and that of simple modules $L(\aalpha)$ 
(or ${}^RL(\aalpha)$) are the same,
by showing that simple modules $L(\aalpha)$ 
are classified according to the equivalence relation $\sim$ in
Theorem \ref{thm:Liffsim}.
In this section, we make another way to prove the coincidence of
the classifications;
we connect $A$-hypergeometric systems $H_A(\aalpha)$
and right $D(R_A)$-modules ${}^R M(\aalpha)$ by functors.
This proof is intrinsic, and hence more desirable.

\subsection{The bimodule $D(R, R_A)$}

In this subsection, we decompose the $(D(R_A), D(R))$-bimodule $D(R, R_A)$
into its $\Z^d$-graded parts similarly to Theorem \ref{thm:Jones}.

Let $\C[x, x^{-1}]$ denote the Laurent polynomial
ring $\C[x_1^{\pm 1},\ldots, x_n^{\pm 1}]$.
By \cite[p. 31]{Cannings},
we have
$$
D(R, R_A)=
\{
P\in D(\C[x, x^{-1}], \C[t, t^{-1}])\,
:\,
P(R)\subseteq R_A
\}.
$$
From \cite[1.3 (e)]{Smith-Stafford},
$$
D(\C[x, x^{-1}], \C[t, t^{-1}])=D(\C[x, x^{-1}])/I_A(x)D(\C[x, x^{-1}]),
$$
and hence
$$
D(\C[x, x^{-1}], \C[t, t^{-1}])=
\bigoplus_{\a\in \Z^d}
t^\a\C[\theta_1,\ldots,\theta_n],
$$
where $\theta_j=x_j\frac{\partial}{\partial x_j}$ ($j=1,\ldots, n$).
From \cite[1.3 (e)]{Smith-Stafford} again,
we have
\begin{equation}
\label{eq:D(R)/I_AD(R)}
D(R, R_A)=D(R)/I_A(x)D(R).
\end{equation}

Note that $D(R_A)\subseteq D(R,R_A)$ by \eqref{eq:D(R)/I_AD(R)}.
Here we identify $t^{\a_j}$ and $s_i$ with $x_j$ and
$\sum_{j=1}^n a_{ij}\theta_j$ respectively
($j=1,\ldots,n$; $i=1,\ldots,d$).
In fact, we have
$$
D(R_A)=D(R, R_A)\cap D(\C[ t^{\pm 1}])
$$
in $D(\C[x, x^{-1}], \C[t, t^{-1}])$.
The bimodule $D(R, R_A)$ inherits the $\Z^d$-grading from 
$D(\C[x, x^{-1}], \C[t, t^{-1}])$,
$$
D(R,R_A)_\a=D(R, R_A)\cap t^\a\C[\theta_1,\ldots,\theta_n].
$$

\begin{proposition}
$$
D(R, R_A)=
\bigoplus_{\a\in\Z^d}
t^\a
\I(\tilde{\Omega}_A(\a)),
$$
where
$\tilde{\Omega}_A(\a)=\N^n\cap T_A^{-1}(\Omega_A(\a))$,
$T_A$ is the linear map from $\Z^n$ to $\Z^d$
defined by $A$,
and $\I(\tilde{\Omega}_A(\a))$ is the ideal of 
$\C[\theta]=\C[\theta_1,\ldots,\theta_n]$ 
vanishing on $\tilde{\Omega}_A(\a)$.
\end{proposition}

\begin{proof}
We have
\begin{eqnarray*}
&&t^\a p(\theta)\in D(R, R_A)_\a\\
&&t^\a p(\theta)(x^\m)\in \C[\N A]\quad (\forall \m\in \N^n)\\
&\Leftrightarrow&p(\m)t^{\a+A\m}\in \C[\N A]\quad (\forall \m\in \N^n)\\
&\Leftrightarrow&
\mbox{$\a+A\m\in \N A$ or $p(\m)=0$ if $\m\in \N^n$}\\
&\Leftrightarrow&
p(\m)=0 \quad (\forall \m\in \N^n\setminus T_A^{-1}(-\a+\N A)).
\end{eqnarray*}
\end{proof}

\begin{corollary}
$D(R, R_A)_\a=t^\a\C[\theta]$
for all $\a\in \N A$.
\end{corollary}

\begin{proof}
In this case $\Omega(\a)=\emptyset$.
Hence $\I(\tilde{\Omega}_A(\a))=\C[\theta]$.
\end{proof}

To close this subsection, we describe the weight decomposition of
${}^R H_A(\aalpha)$.
Let
$$
{}^R H_A(\aalpha)_\llambda:=\{
x\in {}^R H_A(\aalpha)\, :\,
x.(\sum_{j=1}^n a_{ij}\theta_j +\lambda_i)=0\quad (i=1,\ldots, d)\}
$$
for $\llambda={}^t(\lambda_1,\ldots, \lambda_d)$.
Note that
the weight space $D(R, R_A)_\a$ is a $(\C[s], \C[\theta])$-bimodule,
since $D(R_A)_\0=\C[s]$ and $D(R)_\0=\C[\theta]$.
We have
\begin{eqnarray}
{}^R H_A(\aalpha)
&=&\bigoplus_{\a\in \Z^d}
{}^R H_A(\aalpha)_{-\aalpha+\a},\nonumber\\
&=&\bigoplus_{\a\in \Z^d}
D(R,R_A)_\a/(s-\aalpha)D(R,R_A)_\a\nonumber\\
&=&\bigoplus_{\a\in \Z^d}
D(R,R_A)_\a/D(R,R_A)_\a (A\theta +\a -\aalpha)\nonumber\\
&=&\bigoplus_{\a\in \Z^d}
t^\a\left(
\I(\widetilde{\Omega}(\a))/
\I(\widetilde{\Omega}(\a))(A\theta +\a -\aalpha)
\right).
\label{eq:WeightH}
\end{eqnarray}

\subsection{Functors}

Let ${\rm Mod}^R(D(R_A))$ denote the category of right $D(R_A)$-modules,
and ${\rm Mod}_A^R(D(R))$ the category of right $D(R)$-modules
supported by the affine toric variety $\V(I_A(x))$ defined by $A$.
A right $D(R)$-module $N$ is said to be supported by $\V(I_A(x))$
if for every $x\in N$ there exists $m\in \N$ such that
$x I_A(x)^m=0$.

Let $\Phi$ denote the functor from ${\rm Mod}^R(D(R_A))$
to ${\rm Mod}_A^R(D(R))$ defined by
$$
\Phi(M):= M\otimes_{D(R_A)}D(R, R_A),
$$
and $\Psi$ the functor from ${\rm Mod}_A^R(D(R))$
to ${\rm Mod}^R(D(R_A))$ defined by
\begin{eqnarray*}
\Psi(N)&:=& {\rm Hom}_{D(R)}( D(R, R_A), N)\\
&=&\{ x\in N\, :\,
xI_A(x)=0\}.
\end{eqnarray*}
Then $\Psi$ is right adjoint to $\Phi$,
$$
{\rm Hom}_{D(R)}(\Phi(M),N)
\simeq
{\rm Hom}_{D(R_A)}(M, \Psi(N)).
$$
\begin{remark}
The functor $\Phi$ is the direct image functor, for right $D$-modules,
of the closed inclusion
$$
\V(I_A(x))=\{ x\in \C^n\,:\,
f(x)=0\quad (\forall f\in I_A(x))\}\to\C^n.
$$
For a closed embedding between nonsingular varieties,
such a direct image functor gives an equivalence of categories,
known as Kashiwara's equivalence (see e.g. 
\cite[Theorem 4.30]{Kashiwara-book}).
In our case, the affine toric variety $\V(I_A(x))$
is singular whenever $n\neq d$, and the cone $\R_{\geq 0}A$ 
generated by $A$ is strongly convex.
\end{remark}
We have
$$
\Phi(D(R_A))=D(R,R_A),
$$
and
$$
\Psi(D(R,R_A))=\End_{D(R)}(D(R)/I_A(x)D(R))=D(R_A).
$$

The following proposition is immediate from the definitions.

\begin{proposition}
\label{prop:MimpliesH}
We have
$$
\Phi({}^R M(\aalpha))={}^R H_A(\aalpha).
$$
Hence, if ${}^R M(\aalpha)\simeq {}^R M(\bbeta)$,
then ${}^R H_A(\aalpha)\simeq {}^R H_A(\bbeta)$.
\end{proposition}

To show the inverse of Proposition \ref{prop:MimpliesH},
We need the following lemma.

\begin{lemma}
\label{lemma:EndH}
$$
\End_{D(R)}({}^R H_A(\aalpha))=\C\, id.
$$
\end{lemma}

\begin{proof}
Let $\psi\in \End_{D(R)}({}^R H_A(\aalpha))$.
Since $\psi(\overline{1})\in {}^R H_A(\aalpha)_{-\aalpha}$,
there exists a polynomial $f(\theta)\in \C[\theta]$
such that $\psi(\overline{1})=\overline{f(\theta)}$.
Here $\overline{P}$ is the element of ${}^R H_A(\aalpha)$
represented by $P\in D(R)$.
Let $\u, \v\in \N^n$ satisfy $A\u=A\v$.
Then
\begin{eqnarray*}
&&0=\psi(\overline{x^\u-x^\v})
=\psi(\overline{1})(x^\u-x^\v)\\
&&\qquad
=\overline{f(\theta)}(x^\u-x^\v)
=\overline{t^{A\u}(f(\theta+\u)-f(\theta +\v))}.
\end{eqnarray*}
By \eqref{eq:WeightH},
we have
$$
f(\theta+\u)-f(\theta+\v)\in (A\theta+A\u -\aalpha)\C[\theta]
$$
for $\u, \v\in \N^n$ with $A\u=A\v$.
Hence $f(\theta)-f(\theta+\l)\in (A\theta-\aalpha)\C[\theta]$
for all $\l\in L$,
where
$L=\{ \l\in \Z^n\, :\, A\l=\0 \}$.
Letting $A\ggamma=\aalpha$ ($\ggamma\in \C^n$),
we have $f(\ggamma)-f(\ggamma+\l)=0$ for all $\l\in L$.
Thus $f(\ggamma+\theta)\in f(\ggamma)+(A\theta)\C[\theta]$,
or equivalently
$$
f(\theta)\in f(\ggamma)+(A\theta-\aalpha)\C[\theta].
$$
Hence $\psi(\overline{1})=\overline{f(\theta)}
=f(\ggamma)$.
Therefore $\psi=f(\ggamma)id$.
\end{proof}

Let $\bbeta-\aalpha\in \Z^d$, and $Q\in D(R_A)_{\bbeta-\aalpha}$. Then
\begin{eqnarray}
&&\phi_Q: {}^R M(\aalpha)\ni \overline{P} \mapsto
\overline{QP} \in {}^R M(\bbeta)\\
&&\psi_Q: {}^R H_A(\aalpha)\ni \overline{P} \mapsto
\overline{QP} \in {}^R H_A(\bbeta)
\end{eqnarray}
are well-defined morphisms.
Clearly $\psi_Q=\Phi(\phi_Q)$.

\begin{lemma}
The natural map
$$
D(R_A)_{\bbeta-\aalpha}\ni Q\mapsto \phi_Q\in
\Hom_{D(R_A)}({}^R M(\aalpha), {}^R M(\bbeta))
$$
is surjective.
\end{lemma}

\begin{proof}
Let $\phi\in
\Hom_{D(R_A)}({}^R M(\aalpha), {}^R M(\bbeta))$.
Since $\phi(\overline{1})\in {}^R M(\bbeta)_{-\aalpha}$,
there exists $Q\in D(R_A)_{\bbeta-\aalpha}$ such that
$\phi(\overline{1})=\overline{Q}$.
Then $\phi=\phi_Q$.
\end{proof}

As to $\Hom_{D(R)}({}^R H_A(\aalpha), {}^R H_A(\bbeta))$,
we have the following by Lemma \ref{lemma:EndH}.

\begin{corollary}
\label{cor:HomH}
Suppose that ${}^R H_A(\aalpha)\simeq {}^R H_A(\bbeta)$.
Then
$$
\dim_\C \Hom_{D(R)}({}^R H_A(\aalpha), {}^R H_A(\bbeta))=1,
$$
and
the natural map
$D(R_A)_{\bbeta-\aalpha}\to 
\Hom_{D(R)}({}^R H_A(\aalpha), {}^R H_A(\bbeta))$
is surjective.
\end{corollary}

\begin{proof}
The first statement is immediate from Lemma \ref{lemma:EndH}.
The second follows from the fact that in this case the image of the map
is not zero by \cite{IsoClass}.
\end{proof}

Now we are in position of proving the inverse of Proposition
\ref{prop:MimpliesH}.
Note that
there exists a natural morphism $M\to \Psi(\Phi(M))$
for $M\in {\rm Mod}^R(D(R_A))$.

\begin{proposition}
\label{prop:HomHM}
Suppose that ${}^R H_A(\aalpha)\simeq {}^R H_A(\bbeta)$. Then
$$
\Hom_{D(R_A)}({}^R M(\aalpha), {}^R M(\bbeta))
\simeq
\Hom_{D(R)}({}^R H_A(\aalpha), {}^R H_A(\bbeta)).
$$
\end{proposition}

\begin{proof}
Corollary \ref{cor:HomH} states that
there exist $Q, R\in D(R_A)$ such that
\begin{eqnarray*}
\psi_Q &:& {}^R H_A(\aalpha)\ni 
\overline{P}\mapsto
\overline{QP} \in {}^R H_A(\bbeta)\\
\psi_R &:& {}^R H_A(\bbeta)\ni 
\overline{P}\mapsto
\overline{RP} \in {}^R H_A(\aalpha)
\end{eqnarray*}
satisfy
$\psi_R\circ \psi_Q = id_{{}^R H_A(\aalpha)}$ and
$\psi_Q\circ \psi_R = id_{{}^R H_A(\bbeta)}$.

The image of
the natural morphism ${}^RM(\aalpha)\to \Psi(\Phi({}^RM(\aalpha)))
=\Psi({}^R H_A(\aalpha))$ is
$1_\aalpha D(R_A)$,
where
$$
\Psi({}^R H_A(\aalpha))=
\Hom_{D(R)}(D(R, R_A), {}^R H_A(\aalpha)) \ni
1_\aalpha:
\overline{P} \mapsto \overline{P}.
$$
The isomorphism $\psi_Q$ induces an isomorphism $\Psi(\psi_Q)$ : 
$$
\Hom_{D(R)}(D(R, R_A), {}^R H_A(\aalpha))\to
\Hom_{D(R)}(D(R, R_A), {}^R H_A(\bbeta)).
$$
We show that the restriction of $\Psi(\psi_Q)$ to $1_\aalpha D(R_A)$
gives an isomorphism $1_\aalpha D(R_A) \simeq 1_\bbeta D(R_A)$.
By definition $\Psi(\psi_Q)(1_\aalpha)=1_\bbeta Q\in 1_\bbeta D(R_A)$.
Hence $\Psi(\psi_Q)(1_\aalpha D(R_A))\subseteq 1_\bbeta D(R_A)$.
In addition, we have
$$
1_\bbeta=1_\bbeta QR=\Psi(\psi_Q)(1_\aalpha R).
$$
Hence $\Psi(\psi_Q)(1_\aalpha D(R_A))= 1_\bbeta D(R_A)$.
We have thus proved that $\psi_Q$ induces an isomorphism
$1_\aalpha D(R_A)\simeq 1_\bbeta D(R_A)$.
It is lifted to an isomorphism between their projective covers
in ${}^R\cO$,
${}^R M(\aalpha)\simeq {}^R M(\bbeta)$,
which is clearly $\phi_Q$.
Thus $\phi_Q$ and
$\psi_Q$ correspond to each other.
\end{proof}

Combining Propositions
\ref{prop:MimpliesH} and \ref{prop:HomHM}, we have the following.

\begin{corollary}
\label{cor:MRiffHR}
$
{}^R M(\aalpha)\simeq {}^R M(\bbeta)
$
if and only if 
$
{}^R H_A(\aalpha)\simeq {}^R H_A(\bbeta).
$
\end{corollary}

We summarize the classification.
\begin{theorem}
\label{thm:summary}
The following are equivalent:
\begin{enumerate}
\item
$\aalpha\sim\bbeta$.
\item
$
H_A(\aalpha)\simeq H_A(\bbeta).
$
\item
$
{}^R H_A(\aalpha)\simeq {}^R H_A(\bbeta).
$
\item
$
M(\aalpha)\simeq M(\bbeta).
$
\item
$
{}^R M(\aalpha)\simeq {}^R M(\bbeta).
$
\item
$
L(\aalpha)\simeq L(\bbeta).
$
\item
$
{}^R L(\aalpha)\simeq {}^R L(\bbeta).
$
\end{enumerate}
\end{theorem}


\section{Primitive ideals}

In general, a left primitive ideal is not necessarily
a right primitive ideal.
In our case, however, we have the following proposition
thanks to the duality.

\begin{proposition}
For each $\aalpha\in \C^d$, the annihilators
$\Ann(L(\aalpha))$ and 
$\Ann({}^R L(\aalpha))$ coincide, i.e.,
$$
\Ann(L(\aalpha))=
\Ann({}^R L(\aalpha)).
$$
\end{proposition}

\begin{proof}
By definition, it is clear that
$$
\Ann\, M\subseteq \Ann\, M^*
$$
for $M\in \cO\, ({}^R\cO)$.
Since ${}^RL(\aalpha)^*= L(\aalpha)$
and $L(\aalpha)^*= {}^R L(\aalpha)$,
the statement follows.

\end{proof}

Hence the set of annihilators of $\Z^d$-graded simple left 
$D(R_A)$-modules and that of $\Z^d$-graded simple right 
$D(R_A)$-modules are the same.
We denote it by
${\rm Prim}(D(R_A))$,
$$
{\rm Prim}(D(R_A))=
\{\, {\rm Ann}(L(\aalpha))\, :\, \aalpha\in\C^d\,\}.
$$

Next we describe the graded components of the annihilator
ideal $\Ann\, L(\aalpha)$.
For $\aalpha\in \C^d$ and $\a\in \Z^d$, we
define a subset $\Lambda_{[\aalpha]}(\a)$ of $\aalpha+\Z^d$ by
$$
\Lambda_{[\aalpha]}(\a)
:=\{\, \bbeta\in \aalpha+\Z^d\, :\, \bbeta\sim\aalpha,\, \bbeta+\a\sim\aalpha\,\}.
$$

\begin{proposition}[Proposition 3.2.2 in \cite{Musson-Van den Bergh}]
\label{theorem:main}

Let $$(\Ann\, L(\aalpha))_\a=\Ann\, L(\aalpha)\cap D(R_A)_\a$$ 
for $\a\in \Z^d$.
Then
$$
(\Ann\, L(\aalpha))_\a=t^\a\I(\Omega(\a)\cup \Lambda_{[\aalpha]}(\a))
$$
for all $\a$.
\end{proposition}

\begin{proof}
Let $\bbeta\sim\aalpha$, and let
$0\neq v_\bbeta\in L(\aalpha)_\bbeta$.
Then
$$
\Ann(v_\bbeta)_\a
=
\left\{
\begin{array}{ll}
D(R_A)_\a & (\bbeta+\a\not\sim\aalpha)\\
D(R_A)_\a\cap t^\a(s-\bbeta)
&
(\bbeta+\a\sim\aalpha)
\end{array}
\right.
$$
Hence
$$
(\Ann\, L(\aalpha))_\a
=\bigcap_{\bbeta\sim\aalpha,\, \bbeta+\a\sim\aalpha}
D(R_A)_\a\cap t^\a(s-\bbeta).
$$
We have thus proved the assertion.
\end{proof}



\section{Finiteness}

In this section, we prove that the set ${\rm Prim}(D(R_A))$ is finite.
If $\aalpha-\bbeta\notin \Z^d$, then $\aalpha\not\sim\bbeta$,
and hence there exist infinitely many isomorphism classes 
of simple objects $L(\aalpha)$.
We however show that if we properly perturb a parameter $\aalpha$
then the annihilator ideal $\Ann\, L(\aalpha)$ remains unchanged.

First we recall the primitive integral support
function of a facet (maximal proper face) of the cone $\mathbb{R}_{\geq 0}A$.
We denote by $\mathcal{F}$ the set of facets of the cone
$\mathbb{R}_{\geq 0}A$.
Given $\sigma\in\mathcal{F}$,
we denote by $F_\sigma$ the primitive integral support function
of $\sigma$, i.e., $F_\sigma$ is a uniquely determined
linear form on $\mathbb{R}^d$ satisfying
\begin{enumerate}
\item
$F_\sigma(\mathbb{R}_{\geq 0}A)\geq 0$,
\item
$F_\sigma(\sigma)=0$,
\item
$F_\sigma(\mathbb{Z}^d)=\mathbb{Z}$.
\end{enumerate}

Given $\aalpha\in \C^d$, put
\begin{eqnarray*}
\mathcal{F}(\aalpha):&=&
\{\, \sigma\in \cF\, :\, F_\sigma(\aalpha)\in \Z \,\},\\
V(\aalpha):&=& \bigcap_{\sigma\in \mathcal{F}(\aalpha)}(F_\sigma=0).\\
\end{eqnarray*}

Clearly $\mmu\in V(\aalpha)$ implies 
$\cF(\aalpha)\subseteq \cF(\aalpha+\mmu)$.
Note that
$V(\aalpha)$ may not be a linearization of a face.

\begin{example}
Let
$$A = (\a_1,\a_2,\a_3,\a_4)=
\left(
\begin{array}{cccc} 
1 & 0 & 0 & 1\\
0 & 1 & 0 & 1\\
0 & 0 & 1 & -1
\end{array} \right). $$
There are four facets:
$\sigma_{23}:=\R_{\geq 0}\a_2+\R_{\geq 0}\a_3$,
$\sigma_{13}:=\R_{\geq 0}\a_1+\R_{\geq 0}\a_3$,
$\sigma_{24}:=\R_{\geq 0}\a_2+\R_{\geq 0}\a_4$,
$\sigma_{14}:=\R_{\geq 0}\a_1+\R_{\geq 0}\a_4$.
Their primitive integral support functions are respectively
$F_{\sigma_{23}}(s)=s_1$,
$F_{\sigma_{13}}(s)=s_2$,
$F_{\sigma_{24}}(s)=s_1 +s_3$,
$F_{\sigma_{14}}(s)=s_2+s_3$.
Let $\aalpha={}^t(\sqrt{2}, 0, -\sqrt{2})$.
Then $\cF(\aalpha)=\{ \sigma_{13}, \sigma_{24}\}$,
and $V(\aalpha)=\{ {}^t(x, 0, -x)\, :\, x\in\C\}$,
which is not a linearization of any face of $\R_{\geq 0}A$.
\end{example}

Next we briefly review the finite sets $E_\tau(\aalpha)$
defined in \cite{IsoClass}.
Associated to a parameter vector $\aalpha\in \C^d$ and a face $\tau$
of the cone $\R_{\geq 0}A$,  $E_\tau(\aalpha)$ was defined by
\begin{equation}
\label{def:Etau}
E_\tau(\aalpha)=\{
\llambda\in \C(A\cap\tau)/\Z(A\cap\tau)\, :\,
\aalpha-\llambda\in \N A+\Z(A\cap\tau)\}.
\end{equation}
The set $E_\tau(\aalpha)$ has at most
${}^\sharp [\Q(A\cap\tau)\cap\Z^d : \Z(A\cap\tau)]$
elements.
Indeed, suppose that $\llambda\in \C(A\cap\tau)$ satisfies
$\aalpha-\llambda\in \Z^d$. Then
$E_\tau(\aalpha)$ is a subset of
\begin{equation}
\label{set:PossibleEtau}
(\llambda+\Q(A\cap\tau)\cap\Z^d)/\Z(A\cap\tau)
\end{equation}
(see \cite[Proposition 2.3]{IsoClass}).
For a facet $\sigma$, $E_\sigma(\aalpha)\neq\emptyset$
if and only if $F_\sigma(\aalpha)\in F_\sigma(\N A)$,
and, for faces $\tau\preceq\tau'$,
$E_{\tau'}(\aalpha)=\emptyset$ implies
$E_{\tau}(\aalpha)=\emptyset$
(see \cite[Proposition 2.2]{IsoClass}).

We have already introduced
a partial ordering $\preceq$ 
into the parameter space $\C^d$ in \eqref{eq:PartialOrder1}. 
The following is its original definition 
in \cite[Definition 4.1.1]{Saito-Traves}:
For $\aalpha,\bbeta\in \C^d$, we write 
\begin{equation}
\label{def:Preceq}
\text{$\aalpha\preceq\bbeta$\qquad
if \qquad $E_\tau(\aalpha)\subseteq E_\tau(\bbeta)$
for all faces $\tau$.}
\end{equation}
Hence
$\aalpha\sim\bbeta$ if and only if $E_\tau(\aalpha)= E_\tau(\bbeta)$
for all faces $\tau$.
Note that $\aalpha+\Z^d$ has only finitely many equivalence classes,
since $E_\tau(\aalpha)$ and $E_\tau(\aalpha+\a)$ ($\a\in \Z^d$)
are subsets of the finite set \eqref{set:PossibleEtau}
for $\llambda\in\C(A\cap\tau)$ with $\aalpha-\llambda\in \Z^d$.

\begin{lemma}
\label{lemma:translation}
Suppose that $\mmu\in V(\aalpha)$ and 
$\mathcal{F}(\aalpha+\mmu)=\mathcal{F}(\aalpha)$.
Let $\tau$ be a face of the cone $\R_{\geq 0}A$.
Then $E_\tau(\aalpha+\mmu)=\emptyset$ if and only if
$E_\tau(\aalpha)=\emptyset$.
Moreover, if $E_\tau(\aalpha)\not= \emptyset$, then
\begin{equation}
E_\tau(\aalpha+\mmu)=\mmu + E_\tau(\aalpha).
\end{equation}
\end{lemma}

\begin{proof}
By symmetry, it is sufficient to prove
that if $\llambda\in E_\tau(\aalpha)$ then
$\mmu\in \C(A\cap\tau)$ and
$\llambda+\mmu\in E_\tau(\aalpha+\mmu)$.

Suppose that $\llambda\in E_\tau(\aalpha)$.
Then $F_\sigma(\aalpha)=F_\sigma(\aalpha-\llambda)
\in F_\sigma(\N A)\subseteq \N$
for all facets $\sigma\succeq\tau$.
Hence $\C(A\cap\tau)$ is the intersection of a subset of $\mathcal{F}(\aalpha)$,
and thus $\C(A\cap\tau)\supseteq V(\aalpha)$.
This proves $\mmu\in \C(A\cap\tau)$.
Since $\aalpha+\mmu -(\llambda+\mmu)=\aalpha-\llambda\in
\N A +\Z(A\cap\tau)$,
we see $\llambda+\mmu\in E_\tau(\aalpha+\mmu)$.
\end{proof}

\begin{corollary}
\label{corollary:translation}
Suppose that $\mmu\in V(\aalpha)$ and 
$\mathcal{F}(\aalpha+\mmu)=\mathcal{F}(\aalpha)$.
Then 
$\aalpha\sim\bbeta$ if and only if
$\aalpha+\mmu \sim \bbeta+\mmu$.
\end{corollary}

\begin{proof}
If $\aalpha-\bbeta\notin \Z A$, then
$\aalpha\not\sim\bbeta$ and $\aalpha+\mmu\not\sim
\bbeta+\mmu$.

Suppose that
$\aalpha-\bbeta\in \Z A$.
Then
$\mathcal{F}(\aalpha)=\mathcal{F}(\bbeta)$,
and
$\mathcal{F}(\aalpha+\mmu)=\mathcal{F}(\bbeta+\mmu)$.
The assertion follows from Lemma 
\ref{lemma:translation}.
\end{proof}

\begin{proposition}
\label{proposition:translation}
Suppose that $\mmu\in V(\aalpha)$ and 
$\mathcal{F}(\aalpha+\mmu)=\mathcal{F}(\aalpha)$.
Then 
\begin{equation}
{\rm Ann}(L(\aalpha))\, =\,
{\rm Ann}(L(\aalpha+\mmu)).
\end{equation}
\end{proposition}

\begin{proof}
Recall Proposition \ref{theorem:main},
$$
(\Ann\, L(\aalpha))_\a=t^\a\I(\Omega(\a)\cup \Lambda_{[\aalpha]}(\a)),
$$
where
$\Lambda_{[\aalpha]}(\a)
=\{\, \bbeta\, :\, \bbeta\sim \bbeta+\a\sim\aalpha\,\}$.
We prove that
$$
{\rm ZC}(\Lambda_{[\aalpha+\mmu]}(\a))={\rm ZC}(\Lambda_{[\aalpha]}(\a)),
$$
where ${\rm ZC}$ stands for Zariski closure.
Since $\Lambda_{[\aalpha+\mmu]}(\a)=\emptyset$ if and only if
$\Lambda_{[\aalpha]}(\a)=\emptyset$ by Corollary \ref{corollary:translation}, we suppose that they are not empty.
Then again
by Corollary \ref{corollary:translation}
$$
\Lambda_{[\aalpha+\mmu]}(\a)=
\mmu +\Lambda_{[\aalpha]}(\a).
$$
Hence for the proof it suffices to show that
\begin{equation}
V(\aalpha)+{\rm ZC}(
\Lambda_{[\aalpha]}(\a))\, =\,
{\rm ZC}(\Lambda_{[\aalpha]}(\a)).
\end{equation}

Let $\bbeta\in \Lambda_{[\aalpha]}(\a)$, and
$\mmu'\in V(\aalpha)\cap\Z A$.
Put
$v:=\prod_{\C\tau\supseteq V(\aalpha)}[\Z A\cap \Q\tau : \Z(A\cap\tau)]$.
We show that $\bbeta+v \mmu'\in \Lambda_{[\aalpha]}(\a)$.

Suppose that $\C\tau\not\supseteq V(\aalpha)$.
Then there exists a facet $\sigma\succeq\tau$ such that
$\C\sigma\not\supseteq V(\aalpha)$.
Thus $F_\sigma(\aalpha)\notin \Z$.
Note that $v\mmu', \a, \aalpha-\bbeta\in \Z A$, since $\aalpha\sim\bbeta$.
Hence $F_\sigma(\bbeta), F_\sigma(\bbeta+v\mmu'), 
F_\sigma(\bbeta +\a +v\mmu')\notin \Z$.
This implies
$E_\sigma(\bbeta)=E_\sigma(\bbeta+v\mmu')=E_\sigma(\bbeta+\a+v\mmu')=
\emptyset = E_\sigma(\aalpha)$, and
$E_\tau(\bbeta)=E_\tau(\bbeta+v\mmu')=E_\tau(\bbeta+\a+v\mmu')=
\emptyset = E_\tau(\aalpha)$ by \cite[Proposition 2.2]{IsoClass}.

Next suppose that $\C\tau \supseteq V(\aalpha)$.
Then
$E_\tau(\bbeta+v\mmu')=E_\tau(\bbeta)$, and
$E_\tau(\bbeta+\a+v\mmu')=
E_\tau(\bbeta +\a)$,
since $v\mmu'\in \Z(A\cap\tau)$.

Hence $\bbeta+v \mmu'\in \Lambda_{[\aalpha]}(\a)$.
We have thus proved that
$$
\Lambda_{[\aalpha]}(\a) +v (V(\aalpha)\cap\Z A)
\subseteq
\Lambda_{[\aalpha]}(\a).
$$
Taking the Zariski closures, we see that
$$
{\rm ZC}(\Lambda_{[\aalpha]}(\a)) + V(\aalpha)
\subseteq
{\rm ZC}(\Lambda_{[\aalpha]}(\a)).
$$
The other inclusion is trivial.
\end{proof}

\begin{lemma}
\label{representative:finite}
Let $\mathcal{F}'$ be a subset of $\mathcal{F}$,
and let $V(\mathcal{F}')$ be the intersection
$\bigcap_{\sigma\in \mathcal{F}'}(F_\sigma=0)$.
Then 
\begin{equation}
\{\, \aalpha\, :\, \mathcal{F}(\aalpha)\supseteq\mathcal{F}'\,\}/
(\Z A + V(\mathcal{F}'))
\end{equation}
is finite.
\end{lemma}

\begin{proof}
The $\Z$-module $\Z A+V(\mathcal{F}')/V(\mathcal{F}')$
is a $\Z$-submodule of
$\{\, \aalpha\, :\, \mathcal{F}(\aalpha)\supseteq\mathcal{F}'\,\}/
V(\mathcal{F}')$.
Both of them are
free of rank $\dim \C^d/V(\mathcal{F}')$.
Hence the index is finite.
\end{proof}

\begin{theorem}
\label{thm:PrimFinite}
The set ${\rm Prim}(D(R_A))$ is finite.
\end{theorem}

\begin{proof}
It suffices to show that
the set
$$
\Prim_{\cF'}
:=\{ \Ann\, L(\aalpha)\, :\, \cF(\aalpha)=\cF'\}
$$
is finite for each $\cF'\subseteq\cF$.

Let $\cF(\aalpha)=\cF'$.
Let $\aalpha_1,\ldots,\aalpha_k$ be a complete set of representatives
of $\{\, \aalpha\, :\, \mathcal{F}(\aalpha)\supseteq\mathcal{F}'\,\}/
(\Z A + V(\mathcal{F}'))$.
We take the representatives so 
that if a coset $\aalpha_j+\Z A+V(\cF')$ has an
element $\bbeta$ with $\cF(\bbeta)=\cF'$ then $\cF(\aalpha_j)=\cF'$.
Then there exist $j$, $\a\in \Z A$, and $\mmu\in V(\cF')$
such that $\aalpha=\aalpha_j+\a+\mmu$.
We have $\cF'=\cF(\aalpha_j)=\cF(\aalpha_j+\a)$.
By Proposition \ref{proposition:translation},
$\Ann\, L(\aalpha)=\Ann\, L(\aalpha_j+\a)$.
Since each $\aalpha_j+\Z A$ has only finitely many 
equivalence classes,
$\Prim_{\cF'}$ is finite.
\end{proof}

We exhibit a computation of $\Prim(D(R_A))$ in Example \ref{ex:11-1}.

\section{Simplicity}

In this section, we discuss the simplicity of $D(R_A)$.
In the first subsection, 
we consider the conditions: the scoredness and the Serre's $(S_2)$.
We prove that the simplicity of $D(R_A)$ implies the scoredness, and that
all $D(R_A)_\a$ are singly generated $\C[s]$-modules if and only if 
the condition $(S_2)$ is satisfied.
In the second subsection, 
we give a necessary and sufficient condition for the simplicity
(Theorem \ref{simple<=>scored,C2}).

We start this section by noting the $\Z^d$-graded version of a well-known fact.

\begin{lemma}
\label{lem:SimplicityAndPrimitiveIdeal}
The ring $D(R_A)$ is simple if and only if
${\Ann\,} L(\aalpha)=\{ 0\}$ for all $\aalpha\in \C^d$.
\end{lemma}

\begin{proof}
First note that
any two-sided ideal of $D(R_A)$ is $\Z^d$-homogeneous.
(An ideal $I$ is said to be $\Z^d$-homogeneous if 
$I=\bigoplus_{\a\in \Z^d}I\cap D(R_A)_\a$.)

It is enough to show that any maximal ideal of $D(R_A)$
is the annihilator of a simple $\Z^d$-graded module.
Let $I$ be a maximal ideal of $D(R_A)$.
Let $J$ be a maximal $\Z^d$-homogeneous left ideal
containing $I$.
Then $D(R_A)/J$ is a simple $\Z^d$-graded $D(R_A)$-module,
and $\Ann(D(R_A)/J)$ contains $I$.
Since $I$ is maximal, we obtain $I=\Ann(D(R_A)/J)$.
\end{proof}

\subsection{Scored semigroups}

We recall the definition of a scored semigroup
(\cite{Saito-Traves}).
The semigroup $\mathbb{N} A$ is said to be {\it scored\,} 
if 
\begin{equation*}
\mathbb{N} A=\bigcap_{\sigma\in\mathcal{F}}
\{\, {\a}\in \mathbb{Z}^d\, :\,
F_\sigma({\a})\in F_\sigma(\mathbb{N} A)\,\}.
\end{equation*}

We know that
$E_\sigma(\a)\neq \emptyset$ if and only if
$F_\sigma(\a)\in F_\sigma(\N A)$ (\cite[Proposition 2.2]{IsoClass}).
Hence
a semigroup $\N A$ is scored if and only if
\begin{equation*}
\N A =
\{\, \a\in\Z^d\, :
\, E_\sigma(\a)\neq \emptyset\quad
\mbox{for all $\sigma\in\mathcal{F}$}\,\}.
\end{equation*}

In the following lemma, we characterize the subset $\N A$ of $\Z^d$
in terms of the finite sets $E_\tau(\a)$.

\begin{lemma}
\label{NA:Etau}
$$
\N A=
\{ \a\in \Z^d\, :\,
{\0}\in E_\tau (\a)\quad
\text{for all faces $\tau$}\}.
$$
\end{lemma}

\begin{proof}
Let $\tau_0$ be the minimal face of $\R_{\geq 0}A$.
We have
\begin{eqnarray*}
&&\{ \a\in \Z A\, :\,
{\0}\in E_\tau (\a)\quad
\text{for all faces $\tau$}\}\\
&&= \bigcap_\tau (\N A +\Z(A\cap\tau))\\
&&= \N A+\Z (A\cap\tau_0).
\end{eqnarray*}
If $\tau_0=\{ {\0}\}$, or if
the cone $\R_{\ge 0}A$ is strongly convex,
then clearly $\N A+\Z (A\cap\tau_0)=\N A$.
Next suppose $\tau_0\not=\{ {\0}\}$.
Then there exist $c_j\in \Z_{>0}$ such that
${\0}=\sum_{a_j\in A\cap\tau_0}c_ja_j$.
Let $\b\in\Z(A\cap\tau_0)$ and
$\b=\sum d_j\a_j$ with $d_j\in\Z$.
Take $N\in\N$ so that $Nc_j+d_j>0$ for all $j$.
Then $\b=\sum(Nc_j+d_j)\a_j\in\N A$.
\end{proof}

Set
\begin{eqnarray}
\label{eqn:S1}
S_1&:=&\{\, \a\in\Z^d\, :\, E_\sigma(\a)\not= \emptyset
\,\, (\forall\sigma\in\cF)\,\},\\
\label{eqn:S2}
S_2&:=&\{\, \a\in\Z^d\, :\, E_\sigma(\a)\ni\0
\,\, (\forall\sigma\in\cF)\,\}.
\end{eqnarray}
Then $S_2=\bigcap_{\sigma\in\cF}(\N A +\Z(A\cap\sigma))$, and
$$
\N A \subset S_2 \subset S_1
$$
by \eqref{eqn:S1}, \eqref{eqn:S2}, and Lemma \ref{NA:Etau}.

\begin{remark}
\begin{enumerate}
\item
Serre's condition $(S_2)$ is the equality
$\N A= S_2$.
\item
The semigroup $\N A$ is scored if and only if $\N A=S_1$.
\item
By the proof of Lemma \ref{NA:Etau}, we have
$$
\N A =
\{\, \a\in\Z^d\, :\, E_{\tau_0}(\a)\ni\0\,\},
$$
where $\tau_0$ is the minimal face.
\end{enumerate}
\end{remark}

Our first aim in this subsection is to show that
the simplicity of $D(R_A)$ implies the scoredness of $\N A$.
We use the following lemma for the proof.

\begin{lemma}
\label{lemma:DimOmega<d}
$$
\mbox{$\dim \ZC(\Omega(\a)) <d$ for all $\a\in \Z^d$.}
$$
\end{lemma}

\begin{proof}
Take $M$ so that 
$$
\{\,\a\in\Z A\, :\, F_\sigma(\a)\ge M\,\}\subseteq \N A
$$
(see e.g. \cite[Lemma 3.6]{Saito-Traves2}).
Then
\begin{equation*}
\ZC\left(
\Omega(\a)
\right)
\subseteq
\bigcup_{\sigma\in\cF,\, F_\sigma(\a)<M}
\bigcup_{m=0}^{M-F_\sigma(\a)-1}(F_\sigma = m).
\end{equation*}
\end{proof}

\begin{proposition}
If $D(R_A)$ is simple,
then
$\N A$ is scored.
\end{proposition}

\begin{proof}
We know
$$
S_1=\bigcap_{\sigma\in\cF}
\{\,\a\in\Z A\, :\, F_\sigma(\a)\in F_\sigma(\N A)\,\}.
$$
Take $M$ as in the proof of Lemma \ref{lemma:DimOmega<d}.
Then
we have
\begin{equation}
\label{S1-NA}
\ZC\left(
S_1\setminus \N A
\right)
\subseteq
\bigcup_{\sigma\in\cF}
\bigcup_{m=0}^{M-1}(F_\sigma = m).
\end{equation}

Suppose that $\N A$ is not scored, and $\aalpha\in S_1\setminus \N A$.
Since $S_1\setminus \N A$ is a union of some equivalence classes,
we have $\Lambda_{[\aalpha]}(\b)\subseteq S_1\setminus \N A$
for all $\b\in \Z A=\Z^d$.
Hence by \eqref{S1-NA}
$\dim \ZC(\Lambda_{[\aalpha]}(\b))<d$ for all $\b$.
By Lemma \ref{lemma:DimOmega<d}, we have
$\dim \ZC(\Omega(\b)\cup \Lambda_{[\aalpha]}(\b))<d$ for all $\b$.
Hence ${\Ann}\, L(\aalpha)\not=0$
by Proposition \ref{theorem:main}.
\end{proof}

\begin{remark}
By Van den Bergh \cite[Theorem 6.2.5]{Van den Bergh},
if $D(R_A)$ is simple, then $R_A$ is Cohen-Macaulay.
Example \ref{ex:NotSimple:CMandScored} shows that $\N A$ being scored 
and $R_A$ being Cohen-Macaulay are not enough for the simplicity
of $D(R_A)$.
\end{remark}

Now we prove the fact announced in Remark \ref{rem:A1andA2}.
\begin{proposition}
\label{A2=S2}
The $\C[s]$-modules
$D(R_A)_{\a}$ are singly generated for all ${\a}\in \Z^d$
if and only if
the semigroup $\N A$ satisfies $(S_2)$.
\end{proposition}

\begin{proof}
First we paraphrase the condition $(S_2)$.
In \cite[Proposition 3.4]{Saito-Traves2},
we have shown that there exist $(\b_i, \tau_i)$ ($i=1,\ldots, l$),
where $\b_i\in \R_{\geq 0}A\cap\Z^d$ and $\tau_i$ is a face of the cone $\R_{\geq 0}A$,
such that
\begin{equation}
\label{Def:StandardDocomp}
(\mathbb{R}_{\geq 0}A\cap\mathbb{Z}^d)\setminus\mathbb{N}A
=
\bigcup_{i=1}^l
(\b_i+\mathbb{Z} (A\cap\tau_i))\cap \mathbb{R}_{\geq 0}A.
\end{equation}
We may assume that this decomposition is irredundant.
Then
$\{ \b_i+\mathbb{Z} (A\cap\tau_i): i=1,\ldots, l\}$
is unique.
By \cite[Lemma 3.6]{Saito-Traves2}, 
for $\sigma\in\mathcal{F}$,
$$
\mathbb{N}A +\mathbb{Z}(A\cap\sigma)=
[\mathbb{R}_{\geq 0}A+\mathbb{R}(A\cap\sigma)]\cap \mathbb{Z}^d
\setminus \bigcup_{\tau_i=\sigma}({\b}_i+\mathbb{Z}(A\cap\tau_i)).
$$
Hence we obtain
$$
\bigcap_{\sigma\in\mathcal{F}}(
\mathbb{N}A +\mathbb{Z}(A\cap\sigma))=
\mathbb{R}_{\geq 0}A\cap \mathbb{Z}^d
\setminus 
\bigcup_{\tau_i\in\mathcal{F}}({\b}_i+\mathbb{Z}(A\cap\tau_i)).
$$
This means that
$\mathbb{N}A$ satisfies $(S_2)$
if and only if
each $\tau_i$ appearing in \eqref{Def:StandardDocomp}
is a facet.

Suppose that $\mathbb{N}A$ satisfies ($S_2$).
Then, by the previous paragraph and \cite[Proposition 5.1]{Saito-Traves2},
we have
\begin{eqnarray*}
{\rm ZC}(\Omega({\a}))
&=&
\bigcup_{F_\sigma({\a})< 0}\,\,
\bigcup_{m<-F_\sigma({\a}),\, m\in F_\sigma(\mathbb{N}A)}
F_\sigma^{-1}(m)\\
&&\qquad\qquad
\cup
\bigcup_{{\b}_i-{\a}\in \mathbb{N}A+\mathbb{Z}(A\cap\tau_i)}
F_{\tau_i}^{-1}({\b}_i-{\a}).
\end{eqnarray*}
Hence $\mathbb{I}(\Omega({\a}))$ is singly generated.

Next suppose that $\mathbb{N}A$ does not satisfy ($S_2$).
Then a face of codimension greater than one appears in
the difference \eqref{Def:StandardDocomp}.
Let $\tau_1$ be a face of codimension greater than one, and let
${\b}_1+\mathbb{Z}(A\cap\tau_1)$ appear in the difference.
Then
$$
{\rm ZC}(\Omega({\b}_1))=
\bigcup_{{\b}_i-{\b}_1\in \mathbb{N}A+\mathbb{Z}(A\cap\tau_i)}
({\b}_i-{\b}_1+\mathbb{C}(A\cap\tau_i)).
$$
We show that $\mathbb{C}(A\cap\tau_1)$ is an irreducible component
of ${\rm ZC}(\Omega({\b}_1))$.
Suppose the contrary. Then there exists $i$ such that
\begin{eqnarray}
\label{eqn:bib1in}
\b_i-\b_1&\in& \N A+\Z(A\cap\tau_i)\\
\C(A\cap\tau_1)&\subseteq& \b_i-\b_1+\C(A\cap\tau_i).
\label{eqn:bib1subset}
\end{eqnarray}
The latter equation \eqref{eqn:bib1subset}
means that $\b_i-\b_1\in \C(A\cap\tau_i)$ and $\tau_1\preceq\tau_i$.
Combining with \eqref{eqn:bib1in},
we have $\b_i-\b_1\in \Z(A\cap\tau_i)$.
This contradicts the irredundancy of \eqref{Def:StandardDocomp}.
We have thus proved $\mathbb{C}(A\cap\tau_1)$ is an irreducible component
of ${\rm ZC}(\Omega({\b}_1))$.
Hence the ideal $\mathbb{I}(\Omega({\b}_1))$ is not singly generated.
\end{proof}

In the rest of this subsection,
we consider the simplicity of $R_A$ and $S_1$.
The following lemma is immediate from
the definition of $E_\tau(\0)$.

\begin{lemma}
\label{0}
$$
E_{\tau}(\0)=\{\0\}
\qquad\text{for all faces $\tau$.}
$$
\end{lemma}

\begin{lemma}
\label{SufficientlyLarge}
Let $\a\in\Z^d$.
Then
there exists $\b\in \N A$ such that
$$
{}^\sharp E_\tau(\a+\b)=
[\Q(A\cap\tau)\cap\Z^d:
\Z(A\cap\tau)].
$$
{\rm (}In this situation, we write $E_\tau(\a+\b)=\text{\rm full}$.{\rm )}
\end{lemma}

\begin{proof}
We may assume that $\a\in \N A$.
Let $\llambda\in \Q(A\cap\tau)\cap\Z^d/
\Z(A\cap\tau)$, 
take its representative, and denote it by $\llambda$ again.
Write
$-\llambda=\sum_k d_k\a_k$
with
$d_k\in\Z$.
Then
$-\sum_{d_k<0}d_k\a_k -\llambda=\sum_{d_k\ge 0}d_k\a_k$.
Hence $\a+(-\sum_{d_k<0}d_k\a_k)-\llambda\in \N A$.
Thus $\llambda\in E_\tau(\a + (-\sum_{d_k<0}d_k\a_k))$.
We repeat this argument for each pair
$(\tau, \llambda)$ to prove the assertion.
\end{proof}

\begin{proposition}
\label{Rsimple:C0}
The semigroup algebra $R_A$ is a simple $\Z^d$-graded $D(R_A)$-module
if and only if
\begin{equation}
\text{
$
\Q (A\cap\tau)\cap\Z^d
=
\Z (A\cap\tau)
$
\qquad for all faces $\tau$.}
\label{C0}
\end{equation}
\end{proposition}

\begin{proof}
By Lemma \ref{0},
$R_A=\C[\N A]$ is a simple graded $D(R_A)$-module
if and only if
$E_\tau(\a)=\{\0\}$ for all faces $\tau$ and $\a\in\N A$.
This is equivalent to the condition \eqref{C0}
by Lemma \ref{SufficientlyLarge}.
\end{proof}

\begin{lemma}
\label{lem:ScoredImpliesC0}
If the semigroup $\N A$ is scored,
then it satisfies \eqref{C0}.
\end{lemma}

\begin{proof}
For a facet $\sigma$, take $M_\sigma\in \N$
so that $M_\sigma$ is greater than any number in 
$\N\setminus F_\sigma(\N A)$.

Let $\tau$ be a face, and let $\x\in \Q(A\cap\tau)\cap\Z^d$.
For each facet $\sigma\not\succeq \tau$,
there exists $\a_\sigma\in A\cap\tau\setminus \sigma$.
Take $m_\sigma$ large enough to satisfy
$F_\sigma(\x)+m_\sigma F_\sigma(\a_\sigma)\geq M_\sigma$.
Let $\y=\x+\sum_{\sigma\not\succeq \tau}m_\sigma \a_\sigma$.
Then $\y\in \Q(A\cap\tau)\cap\Z^d$, and
$$
\begin{array}{ll}
F_\sigma(\y)=0 & \text{if $\sigma\succeq\tau$},\\
F_\sigma(\y)\geq M_\sigma & \text{otherwise}.
\end{array}
$$
Since $\N A$ is scored,
$\y\in \N A$.
Hence $\y\in \N(A\cap\tau)$, and $\x\in \Z(A\cap\tau)$.
\end{proof}

\begin{corollary}
\label{Rsimple:scored}
If the semigroup $\N A$ is scored,
then $R_A$ is a simple $\Z^d$-graded $D(R_A)$-module.
\end{corollary}

\begin{proof}
This follows from Proposition \ref{Rsimple:C0}
and Lemma \ref{lem:ScoredImpliesC0}.
\end{proof}

\begin{proposition}
The semigroup $\N A$ is scored if and only if
$\C [S_1]$ is a simple $\Z^d$-graded $D(R_A)$-module.
\end{proposition}

\begin{proof}
Suppose that $\N A$ is scored.
Then $S_1=\N A$. Hence
$\C [S_1]$ is a simple $\Z^d$-graded $D(R_A)$-module by
Corollary \ref{Rsimple:scored}.

Suppose that $\C [S_1]$ is a simple $\Z^d$-graded $D(R_A)$-module.
Since $R_A$ is a nonzero $\Z^d$-graded $D(R_A)$-submodule 
of $\C [S_1]$, we have $R_A=\C[S_1]$.
Hence $\N A$ is scored.
\end{proof}

\subsection{Conditions for simplicity}

The aim of this subsection is to give a necessary and sufficient condition
for the vanishing of a primitive ideal $\Ann\, L(\aalpha)$.
It leads to a necessary and sufficient condition
for the simplicity of $D(R_A)$.

We start this subsection by introducing some notation.
Let $\aalpha\in \C^d$. Set
\begin{eqnarray*}
\cF_+(\aalpha)&:=&
\{\,\sigma\in\cF\, :\, F_\sigma(\aalpha)\in F_\sigma(\N A)\,\},\\
\cF_-(\aalpha)&:=&
\{\,\sigma\in\cF\, :\, 
F_\sigma(\aalpha)\in \Z\setminus F_\sigma(\N A)\,\},
\end{eqnarray*}
and
$$
\R_{>0}(\aalpha):=
\left\{\, \ggamma\in \R^d\,:\,
\begin{array}{l}
F_\sigma(\ggamma)> 0\quad (\sigma\in\cF_+(\aalpha))\\
F_\sigma(\ggamma)< 0
\quad (\sigma\in\cF_-(\aalpha))
\end{array}
\,\right\}.
$$
Let ${\rm Face}(\aalpha)$ denote the set of faces $\tau$
such that $\aalpha-\llambda\in \Z^d$
for some $\llambda\in \C (A\cap\tau)$,
and that every facet $\sigma$
containing $\tau$ belongs to $\cF_+(\aalpha)$.
Let $[\aalpha]$ denote the equivalence class that $\aalpha$ belongs to.
An equivalence class $[\aalpha]$ is said to be {\it extreme\,}
if 
$E_\tau (\aalpha)$ has $[\Q(A\cap\tau)\cap\Z^d:\Z(A\cap\tau)]$-many
elements (i.e., $E_\tau (\aalpha)=\mbox{full}$)
for every $\tau\in {\rm Face}(\aalpha)$,
and that
$E_\tau (\aalpha)$ is empty for every $\tau\notin {\rm Face}(\aalpha)$.

We compare the conditions:
\begin{enumerate}

\item
An equivalence class $[\aalpha]$ is extreme,

\item
$\R_{>0}(\aalpha)$ is not empty,

\item
$\ZC ([\aalpha])=\C^d$,

\item
$\Ann \, L(\aalpha)=0$.
\end{enumerate}

\begin{remark}
The conditions (1) and (2) have an advantage over the condition (3),
for to check (1) and (2) we do not need the equivalence class $[\aalpha]$,
which is not easy to compute.
\end{remark}

We need the following technical lemmas.

\begin{lemma}
\label{FullEtau:Lattice}
Let $\tau$ be a face of $\R_{\geq 0}A$.
Then there exists $M\in\N$ such that,
if $\a\in \Z^d$ satisfies 
$F_\sigma(\a)\geq M$ for all facets $\sigma\succeq \tau$,
then
$
E_\tau(\a)
=
\text{\rm full}.
$
\end{lemma}

\begin{proof}
Let $\llambda\in \Q (A\cap\tau)\cap \Z^d$.
By \cite[Lemma 3.6]{Saito-Traves2},
there exists $M\in \N$ such that
$\c\in \N A +\Z(A\cap\tau)$
for all $\c\in\Z^d$ satisfying
$F_\sigma(\c)\geq M$ for all facets $\sigma\succeq \tau$.
Hence, if $\a\in \Z^d$ satisfies 
$F_\sigma(\a)\geq M$ for all facets $\sigma\succeq \tau$,
then $\a-\llambda\in \N A +\Z(A\cap\tau)$, or
$\llambda\in E_\tau(\a)$.
\end{proof}

\begin{lemma}
\label{FullEtau:General}
Let $\tau$ be a face of $\R_{\geq 0}A$,
and let $\aalpha\in \C^d$.
Assume that
there exists $\llambda\in \C(A\cap\tau)$
such that $\aalpha-\llambda$ belongs to $\Z^d$.
Then there exists $M\in\N$ such that,
if $\ggamma\in\aalpha+\Z^d$ satisfying
$F_\sigma(\ggamma)\in\Z_{\geq M}$ for all facets $\sigma\succeq \tau$,
then
$
E_\tau(\ggamma)
=
\text{\rm full}.
$
\end{lemma}

\begin{proof}
Apply Lemma \ref{FullEtau:Lattice} to $\ggamma-\llambda$.
\end{proof}

\begin{proposition}
\label{FullDimEquivClass:C2}
If the Zariski closure of 
an equivalence class $[\aalpha]$
is the whole space $\C^d$,
then
$\R_{>0}(\aalpha)$ is not empty.
\end{proposition}

\begin{proof}
Since the set
\begin{equation}
\label{eqn:alpha+ZA}
\left\{\, \ggamma\in \aalpha+\Z^d\,:\,
\begin{array}{l}
F_\sigma(\ggamma)\in F_\sigma(\N A)\qquad (\sigma\in\cF_+(\aalpha))\\
F_\sigma(\ggamma)\in \Z\setminus
F_\sigma(\N A)\quad (\sigma\in\cF_-(\aalpha))
\end{array}
\,\right\}
\end{equation}
contains $[\aalpha]$,
its Zariski closure equals $\C^d$.
Take a real number $\epsilon$
so that
$\epsilon$ is algebraically independent over
$\Q[\, F_\sigma({\rm Re}(\aalpha)),\,
F_\sigma({\rm Im}(\aalpha))\,:\, \sigma\in\cF\,]$.
Put
$\aalpha_\epsilon:={\rm Re}(\aalpha)+\epsilon{\rm Im}(\aalpha)$,
where ${\rm Re}(\aalpha)$ and ${\rm Im}(\aalpha)$ are the vectors
in $\R^d$ with $\aalpha={\rm Re}(\aalpha)+\sqrt{-1}{\rm Im}(\aalpha)$
Then
$\cF_+(\aalpha)=\cF_+(\aalpha_\epsilon)$,
and
$\cF_-(\aalpha)=\cF_-(\aalpha_\epsilon)$,
since $\sigma\in \cF_{\pm}(\aalpha)$ or $\cF_{\pm}(\aalpha_\epsilon)$ 
implies $F_\sigma({\rm Im}(\aalpha))=0$.
The set 
\begin{equation}
\label{eqn:alpha_epsilon+ZA}
\left\{\, \ggamma\in \aalpha_\epsilon+\Z^d\,:\,
\begin{array}{l}
F_\sigma(\ggamma)\in F_\sigma(\N A)\qquad (\sigma\in\cF_+(\aalpha))\\
F_\sigma(\ggamma)\in \Z\setminus
F_\sigma(\N A)\quad (\sigma\in\cF_-(\aalpha))
\end{array}
\,\right\}
\end{equation}
is bijective to the set \eqref{eqn:alpha+ZA}
under the map sending $\aalpha_\epsilon+\a$ to $\aalpha +\a$,
and hence its Zariski closure
equals $\C^d$.
The Zariski closure of 
$$
\mbox{\eqref{eqn:alpha_epsilon+ZA}$\setminus
\bigcup_{\sigma\in \cF_+(\aalpha)}(F_\sigma =0)
\setminus
\bigcup_{\sigma\in \cF_-(\aalpha)}
\bigcup_{m\in \N\setminus F_\sigma(\N A)}(F_\sigma =m)$}
$$
also equals $\C^d$.
Hence
$\R_{>0}(\aalpha)$
is not empty.
\end{proof}

\begin{proposition}
\label{FullDimEquivClass:Rigidity}
If the Zariski closure of an equivalence class $[\aalpha]$
is the whole space $\C^d$,
then $[\aalpha]$ is extreme, i.e.,
$$
[\aalpha]=
\left\{\, \ggamma\in\aalpha+\Z^d\, :\,
\begin{array}{l}
E_\tau (\ggamma)= \text{{\rm full}}\quad (\tau\in {\rm Face}(\aalpha))\\
E_\tau (\ggamma)=\emptyset\qquad (\tau\notin {\rm Face}(\aalpha))
\end{array}\,\right\}.
$$
\end{proposition}

\begin{proof}
By Lemma \ref{FullEtau:General} the equivalence class
$$
[\aalpha_1]:=
\left\{\, \ggamma\in\aalpha+\Z^d\, :\,
\begin{array}{l}
E_\tau (\ggamma)= \text{full}\quad (\tau\in {\rm Face}(\aalpha))\\
E_\tau (\ggamma)=\emptyset\qquad (\tau\notin {\rm Face}(\aalpha))
\end{array}\,\right\}
$$
contains 
\begin{equation}
\label{eqn:FullCore}
\left\{\, \ggamma\in\aalpha+\Z^d\, :\,
\begin{array}{l}
F_\sigma (\ggamma)\geq M \quad (\sigma\in \cF_+(\aalpha))\\
F_\sigma (\ggamma)<0 \qquad (\sigma\in \cF_-(\aalpha))
\end{array}\,\right\}
\end{equation}
for some $M$ sufficiently large.

Suppose that $[\aalpha]\not=
[\aalpha_1]$.
Then $[\aalpha]$ does not belong to the set \eqref{eqn:FullCore}.
Hence we have
$$
[\aalpha]
\subseteq
\bigcup_{\sigma\in\cF_+(\aalpha)}
\bigcup_{m=0}^{M-1}
\{\, \ggamma\in\aalpha+\Z A\,:\,
F_\sigma(\ggamma)=m\,\}.
$$
This contradicts the assumption that the dimension of
$\ZC ([\aalpha])$ equals $d$.
\end{proof}

\begin{proposition}[cf. Proposition 3.3.1 in \cite{Musson-Van den Bergh}]
\label{Ann=0:dim=d}
Let $\aalpha\in\C^d$.
Then
${\Ann}(L(\aalpha))=0$ if and only if
$\ZC ([\aalpha])=\C^d$.
\end{proposition}

\begin{proof}
Let $I:={\Ann}(L(\aalpha))$.
Recall that we have
$$I_{\a}=t^\a\I(\Omega(\a)\cup\Lambda_{[\aalpha]}(\a)),$$
where 
$$
\Lambda_{[\aalpha]}(\a)=
\{\, \ggamma\,:\, \ggamma\sim\aalpha,\, \ggamma+\a\sim\aalpha\,\}.
$$

Since $I_{\0}=\I([\aalpha])$, the vanishing of $I$ leads to
the assertion that $\ZC ([\aalpha])=\C^d$.

Next suppose that $\ZC ([\aalpha])=\C^d$.
As in the proof of Proposition \ref{FullDimEquivClass:Rigidity}
there exists $M\in\N$ such that
$$
[\aalpha]\supseteq
\left\{\, \ggamma\in\aalpha+\Z^d\, :\,
\begin{array}{l}
F_\sigma (\ggamma)\geq M \quad (\sigma\in \cF_+(\aalpha))\\
F_\sigma (\ggamma)<0 \qquad (\sigma\in \cF_-(\aalpha))
\end{array}\,\right\}.
$$
Hence
\begin{eqnarray*}
&&\Lambda_{[\aalpha]}(\a)\\
&&\supseteq
\left\{\, \ggamma\in\aalpha+\Z^d\, :\,
\begin{array}{l}
F_\sigma (\ggamma)\geq \max \{ M, -F_\sigma(\a)\}
 \quad (\sigma\in \cF_+(\aalpha))\\
F_\sigma (\ggamma)<\min\{ 0, -F_\sigma(\a)\} \qquad (\sigma\in \cF_-(\aalpha))
\end{array}\,\right\}.
\end{eqnarray*}
Since the right hand side is $d$-dimensional
by Proposition \ref{FullDimEquivClass:C2},
the Zariski closure
$\ZC (\Lambda_{[\aalpha]}(\a))$ is also $d$-dimensional.
Hence $I_\a=0$ for all $\a\in\Z^d$.
\end{proof}

\begin{proposition}
\label{C1-2=>simple}
If $[\aalpha]$ is extreme, and
$\R_{>0}(\aalpha)$ is not empty,
then $\ZC ([\aalpha])=\C^d$.
\end{proposition}

\begin{proof}
As in the proof of Proposition \ref{FullDimEquivClass:Rigidity},
$[\aalpha]$
contains 
$$
\left\{\, \ggamma\in\aalpha+\Z^d\, :\,
\begin{array}{l}
F_\sigma (\ggamma)> M \quad (\sigma\in \cF_+(\aalpha))\\
F_\sigma (\ggamma)<0 \qquad (\sigma\in \cF_-(\aalpha))
\end{array}\,\right\}
$$
for some $M$ sufficiently large.
By the assumption,
the dimension of
$$
\left\{\, \a\in \Z^d\, :\,
\begin{array}{l}
F_\sigma (\a)> M -F_\sigma(\aalpha)
 \quad (\sigma\in \cF_+(\aalpha))\\
F_\sigma (\a)< -F_\sigma(\aalpha) \qquad (\sigma\in \cF_-(\aalpha))
\end{array}\,\right\}.
$$
equals $d$.
Hence the proposition follows.
\end{proof}

\begin{theorem}
\label{thm:VanishingOfPrimitiveIdeal}
Let $\aalpha\in\C^d$.
Then
${\Ann}(L(\aalpha))=0$ if and only if
$[\aalpha]$ is extreme, and $\R_{>0}(\aalpha)$ is not empty.
\end{theorem}

\begin{proof}
This follows from
Propositions 
\ref{FullDimEquivClass:C2}, \ref{FullDimEquivClass:Rigidity},
and
\ref{C1-2=>simple}.
\end{proof}

\begin{theorem}
\label{simple<=>C1,C2}
The algebra $D(R_A)$ is simple if and only if the
conditions 
\begin{enumerate}

\item[{\rm (C1)}]
Any equivalence class is extreme.

\item[{\rm (C2)}]
For any $\aalpha$,
$\R_{>0}(\aalpha)$ is not empty.
\end{enumerate}
are satisfied.
\end{theorem}

\begin{proof}
This follows from
Lemma \ref{lem:SimplicityAndPrimitiveIdeal}
and
Theorem \ref{thm:VanishingOfPrimitiveIdeal}.
\end{proof}

\begin{remark}
To know whether $D(R_A)$ is simple or not,
by
Theorem \ref{thm:PrimFinite},
we need to check (C1) and (C2) only for 
finitely many $\aalpha$.
\end{remark}

\begin{proposition}
\label{scored=>C1}
If the semigroup $\N A$ is scored, then
it satisfies the condition {\rm (C1)}.
\end{proposition}

\begin{proof}
First note that the condition \eqref{C0} is satisfied in the scored case
by Lemma \ref{lem:ScoredImpliesC0}.

Let $\llambda\in\C(A\cap\tau)$ and $\aalpha-\llambda\in\Z^d$.
Suppose that $\sigma\in\cF_+(\aalpha)$
for all facets $\sigma$ containing $\tau$.
We need to show $\llambda\in E_\tau(\aalpha)$.

If a facet $\sigma$ contains $\tau$, then
$\sigma\in \cF_+(\aalpha)$, or $F_\sigma(\aalpha)\in F_\sigma(\N A)$.
Hence $F_\sigma(\aalpha-\llambda)\in F_\sigma(\N A)$.
If a facet $\sigma$ does not contain $\tau$, then
there exists $\a_j\in A\cap\tau$ such that $F_\sigma(\a_j)>0$.
Hence $F_\sigma(\aalpha-\llambda+m\a_j)\in F_\sigma(\N A)$
for $m\in \N$ sufficiently large.
Hence there exists $\a\in \N(A\cap\tau)$ such that
$F_\sigma(\aalpha-\llambda+\a)\in F_\sigma(\N A)$
for all $\sigma\in\cF$.
Since $\N A$ is scored, we obtain
$\aalpha-\llambda+\a\in \N A$.
This means $\llambda\in E_\tau(\aalpha)$.
\end{proof}

\begin{theorem}
\label{simple<=>scored,C2}
The algebra $D(R_A)$ is simple if and only if
the semigroup $\N A$ is scored and satisfies the condition {\rm (C2)}.
\end{theorem}

\begin{proof}
This immediately follows from
Theorem \ref{simple<=>C1,C2} and Proposition \ref{scored=>C1}.
\end{proof}

\begin{corollary}
Assume that the cone $\R_{\geq 0}A$ is simplicial.
Then the algebra
$D(R_A)$ is simple if and only if
$\N A$ is scored.
\end{corollary}

\begin{proof}
In this case the cone $\R_{\geq 0}A$ has exactly $d$ facets.
Since the $d$ $F_\sigma$'s are linearly independent,
the condition (C2) is satisfied.
\end{proof}
\section{Examples}

\begin{example}
\label{ex:11-1}
Let
$$A = 
\left( \begin{array}{llll} 1 & 1 & 2 & 2 \\
 1 & 2 & 0 & 1 \end{array} \right). $$
Then
$\N A$ is the set of black dots in {\sc Figure 1},
and $\R_{\geq 0}A$ has two facets:
$\sigma_2:=\R_{\geq 0}\a_2$ and $\sigma_3:=\R_{\geq 0}\a_3$.
Their primitive integral support functions are
$F_{\sigma_2}(s)=2s_1-s_2$ and
$F_{\sigma_3}(s)=s_2$ respectively.

\begin{figure}[!ht]
\label{fig:1}
\begin{picture}(50,65)(0,20)
\put(10,10){\circle*{5}}
\put(20,10){\circle{5}}
\put(30,10){\circle*{5}}
\put(40,10){\circle{5}}
\put(50,10){\circle*{5}}
\put(20,20){\circle*{5}}
\put(30,20){\circle*{5}}
\put(40,20){\circle*{5}}
\put(50,20){\circle*{5}}
\put(20,30){\circle*{5}}
\put(30,30){\circle*{5}}
\put(40,30){\circle*{5}}
\put(50,30){\circle*{5}}
\put(30,40){\circle*{5}}
\put(40,40){\circle*{5}}
\put(50,40){\circle*{5}}
\put(30,50){\circle*{5}}
\put(40,50){\circle*{5}}
\put(50,50){\circle*{5}}
\put(10,10){\vector(1,0){50}}
\put(10,10){\vector(0,1){50}}
\put(10,10){\vector(1,2){25}}
\put(37,70){$\sigma_2$}
\put(70,8){$\sigma_3$}
\end{picture}
\caption{The semigroup $\N A$}
\end{figure}

The condition (C2) is satisfied, since $\R_{\geq 0}A$ is simplicial.
However $D(R_A)$ is not simple, since $\N A$ is not scored.
We have
\begin{eqnarray*}
\{\aalpha\, :\, \cF(\aalpha)=\emptyset\}/\C^2
&=&
\{ {}^t(\sqrt{2}, \sqrt{3})\}\\
\{\aalpha\, :\, F_{\sigma_2}(\aalpha)\in \Z\}/(\Z^2+(F_{\sigma_2}=0))
&=&
\{ {}^t(\frac{1}{2}\sqrt{2}, \sqrt{2})\}\\
\{\aalpha\, :\, F_{\sigma_3}(\aalpha)\in \Z\}/(\Z^2+(F_{\sigma_3}=0))
&=&
\{ {}^t(\sqrt{2}, 0)\}\\
\{\aalpha\, :\, F_{\sigma_2}(\aalpha), F_{\sigma_3}(\aalpha)\in \Z\}/\Z^2
&=&
\{ {}^t(0, 0), {}^t(\frac{1}{2}, 0)\}.
\end{eqnarray*}

First we classify $\Z^2$. Let $\aalpha\in \Z^2$.
Then we see
\begin{itemize}
\item
$E_{\R_\geq 0}(\aalpha)=\{ \0\}$,
\item
$E_{\sigma_2}(\aalpha)=\{ \0\} \Leftrightarrow 2\alpha_1-\alpha_2\geq 0$,
\item
$E_{\sigma_3}(\aalpha)=\{ \0, {}^t(1,0)\} \Leftrightarrow \alpha_2\geq 1$,
\item
$E_{\sigma_3}(\aalpha)=\{ \0\} \Leftrightarrow \alpha_2=0$, $\alpha_1\in 2\Z$,
\item
$E_{\sigma_3}(\aalpha)=\{ {}^t(1,0)\} \Leftrightarrow \alpha_2=0$, 
$\alpha_1\in 2\Z+1$,
\item
$E_{\{ \0\}}(\aalpha)=\{ \0\} \Leftrightarrow \aalpha\in \N A$.
\end{itemize}

There are eight classes in $\Z^2$:
\begin{enumerate}
\item
\begin{eqnarray*}
&& \{\aalpha\in \Z^2\,:\,
E_{\sigma_2}(\aalpha)=\{ \0\},\, 
E_{\sigma_3}(\aalpha)=\{ \0, {}^t(1,0)\},\,
E_{\{ \0\}}(\aalpha)=\{ \0\}\}\\
&=&\{ {}^t(\alpha_1,\alpha_2)\in \Z^2\,:\,
\alpha_2\geq 1,\, 2\alpha_1-\alpha_2\geq 0\}.
\end{eqnarray*}
\item
\begin{eqnarray*}
&& \{\aalpha\in \Z^2\,:\,
E_{\sigma_2}(\aalpha)=\emptyset,\, 
E_{\sigma_3}(\aalpha)=\{ \0, {}^t(1,0)\},\,
E_{\{ \0\}}(\aalpha)=\emptyset\}\\
&=&\{ {}^t(\alpha_1,\alpha_2)\in \Z^2\,:\,
\alpha_2\geq 1,\, 2\alpha_1-\alpha_2< 0\}.
\end{eqnarray*}
\item
\begin{eqnarray*}
&& \{\aalpha\in \Z^2\,:\,
E_{\sigma_2}(\aalpha)=\{ \0\},\, 
E_{\sigma_3}(\aalpha)=\{ \0\},\, 
E_{\{ \0\}}(\aalpha)=\{ \0\}\}\\
&=&\{ {}^t(\alpha_1,\alpha_2)\in \Z^2\,:\,
\alpha_2=0,\, \alpha_1\in 2\N\}.
\end{eqnarray*}
\item
\begin{eqnarray*}
&& \{\aalpha\in \Z^2\,:\,
E_{\sigma_2}(\aalpha)=\{ \0\},\, 
E_{\sigma_3}(\aalpha)=\{ {}^t(1,0)\},\,
E_{\{ \0\}}(\aalpha)=\emptyset\}\\
&=&\{ {}^t(\alpha_1,\alpha_2)\in \Z^2\,:\,
\alpha_2=0,\, \alpha_1\in 2\N +1\}.
\end{eqnarray*}
\item
\begin{eqnarray*}
&& \{\aalpha\in \Z^2\,:\,
E_{\sigma_2}(\aalpha)=\emptyset,\, 
E_{\sigma_3}(\aalpha)=\{ \0\},\,
E_{\{ \0\}}(\aalpha)=\emptyset\}\\
&=&\{ {}^t(\alpha_1,\alpha_2)\in \Z^2\,:\,
\alpha_2=0,\, \alpha_1\in 2(-\N-1)\}.
\end{eqnarray*}
\item
\begin{eqnarray*}
&& \{\aalpha\in \Z^2\,:\,
E_{\sigma_2}(\aalpha)=\emptyset,\, 
E_{\sigma_3}(\aalpha)=\{ {}^t(1,0)\},\,
E_{\{ \0\}}(\aalpha)=\emptyset.\}\\
&=&\{ {}^t(\alpha_1,\alpha_2)\in \Z^2\,:\,
\alpha_2=0,\, \alpha_1\in -2\N-1\}.
\end{eqnarray*}
\item
\begin{eqnarray*}
&& \{\aalpha\in \Z^2\,:\,
E_{\sigma_2}(\aalpha)=\{ \0\},\, 
E_{\sigma_3}(\aalpha)=\emptyset,\,
E_{\{ \0\}}(\aalpha)=\emptyset\}\\
&=&\{ {}^t(\alpha_1,\alpha_2)\in \Z^2\,:\,
\alpha_2<0,\, 2\alpha_1-\alpha_2\geq 0\}.
\end{eqnarray*}
\item
\begin{eqnarray*}
&& \{\aalpha\in \Z^2\,:\,
E_{\sigma_2}(\aalpha)=\emptyset,\, 
E_{\sigma_3}(\aalpha)=\emptyset,\,
E_{\{ \0\}}(\aalpha)=\emptyset\}\\
&=&\{ {}^t(\alpha_1,\alpha_2)\in \Z^2\,:\,
\alpha_2<0,\, 2\alpha_1-\alpha_2<0\}.
\end{eqnarray*}
\end{enumerate}

Let $\aalpha$ be not extreme, i.e.,
let $\aalpha$ belong to (3), (4), (5), or (6).
Then
$\ZC(\Lambda_{[\aalpha ]}(\a))=\{\,\mmu\, :\, \mu_2=0\,\}$ if $a_2=0$ and
$a_1\in 2\Z$,
and $\Lambda_{[\aalpha ]}(\a)=\emptyset$ otherwise.

Let $\aalpha$ be extreme, i.e.,
let $\aalpha$ belong to (1), (2), (7), or (8).
Then
$\ZC(\Lambda_{[\aalpha ]}(\a))=\C^2$ for all $\a\in \Z^2$.

Similarly
${}^t(\sqrt{2}, \sqrt{3})+\Z^2$,
${}^t(\frac{1}{2}\sqrt{2}, \sqrt{2})+\Z^2$,
${}^t(\sqrt{2}, 0)+\Z^2$, and
${}^t(\frac{1}{2}, 0)+\Z^2$
have one, two, four, and eight equivalence classes respectively.
We see that $\aalpha={}^t(\alpha_1,\alpha_2)$ is not extreme
if and only if $\alpha_2=0$ if and only if $\Ann\, L(\aalpha)=
\Ann \, L(\0)$.

Hence
$$
{\rm Prim}(D(R_A))=
\{\, (0),\, \Ann\, L(\0) \,\}.
$$
\end{example}

\begin{example}
\label{ex:NotSimple:CMandScored}
(cf. \cite[Example 4.9]{Trung-Hoa}.)
Let
$$A = 
\left( \begin{array}{cccccc} 
1 & 1 & 1 & 1 & 1 & 1\\
0 & 2 & 3 & 0 & 2 & 3\\ 
0 & 0 & 0 & 1 & 1 & 1 
\end{array} \right). 
$$
Then the cone $\R_{\geq 0}A$ has four facets whose primitive
integral support functions are
$$
\begin{array}{ll}
F_{\sigma_{14}}(s)=s_2,\qquad & F_{\sigma_{36}}(s)=3s_1-s_2,\\
F_{\sigma_{123}}(s)=s_3, & F_{\sigma_{456}}(s)=s_1-s_3.
\end{array}
$$
We have an equality
\begin{equation}
\label{eqn:RelationFs}
F_{\sigma_{14}}+F_{\sigma_{36}}
=3(F_{\sigma_{123}}+F_{\sigma_{456}}).
\end{equation}
Then the semigroup
$$
S:=\N A=
\{ \a\in \R_{\geq 0}A\cap\Z^3\, :\, F_{\sigma_{14}}(\a)\neq 1\}
$$
is scored.
Let $\aalpha={}^t(0,1,0)$.
Then ${\rm Face} (\aalpha)=\{ F_{\sigma_{123}}, F_{\sigma_{456}}\}$.
We see that $\aalpha$ is extreme, and $\R_{>0}(\aalpha)$ is empty.
Hence $\Ann \, L(\aalpha)\neq 0$.

Let $\lambda:=F_{\sigma_{14}}+2F_{\sigma_{36}}=6s_1-s_2$,
and let $E_\lambda$ be the blow-up extension of $S$,
$$
E_\lambda:=\{\, 
\begin{pmatrix}
\a\\ 
p
\end{pmatrix}
\in S\bigoplus\N\, :\, 
p\leq \lambda(\a)\,\}.
$$
We can prove that an affine semigroup is scored 
if and only if its blow-up extension is scored by the same argument as the proof of \cite[Lemma 1.1]{Trung-Hoa}.
Thus $E_\lambda$ is scored.
Indeed, the cone $\R_{\geq 0}E_\lambda$ has six facets
whose primitive integral support functions are
$$
\begin{array}{ll}
F_{\widetilde{\sigma_{14}}}(s,p)=F_{\sigma_{14}}(s)  & 
F_{\widetilde{\sigma_{36}}}(s,p)=F_{\sigma_{36}}(s),\\
F_{\widetilde{\sigma_{123}}}(s,p)=F_{\sigma_{123}}(s), & 
F_{\widetilde{\sigma_{456}}}(s,p)=F_{\sigma_{456}}(s),\\
F_{\sigma_p}(s,p)=p & F_{\sigma_\lambda}(s,p)=\lambda(s)-p,
\end{array}
$$
and
$$
E_\lambda=\R_{\geq 0}E_\lambda\cap\Z^4\setminus (F_{\widetilde{\sigma_{14}}}=1).
$$
In addition, $\C[E_\lambda]$ is Cohen-Macaulay \cite[Example 4.9]{Trung-Hoa}.

However $D(\C[E_\lambda])$ is not simple. To see this, let
$\bbeta={}^t(0,1,0,0)$.
Then
$$
\cF_+(\bbeta)=\{ \widetilde{\sigma_{123}}, \widetilde{\sigma_{456}},
\sigma_p\},\quad 
\cF_-(\bbeta)=\{ \widetilde{\sigma_{14}}, \widetilde{\sigma_{36}},
\sigma_\lambda\}.
$$
By \eqref{eqn:RelationFs}
$$
(F_{\widetilde{\sigma_{123}}}>0)\cap(F_{\widetilde{\sigma_{456}}}>0)
\cap(F_{\widetilde{\sigma_{14}}}<0)\cap(F_{\widetilde{\sigma_{36}}}<0)
$$
is empty.
Hence $\R_{>0}(\bbeta)$ is also empty, and
$D(\C[E_\lambda])$ is not simple by 
Theorem \ref{simple<=>scored,C2}.
\end{example}


\bigskip

Department of Mathematics

Hokkaido University

Sapporo, 060-0810

Japan

e-mail: saito@math.sci.hokudai.ac.jp

\end{document}